\documentclass[a4paper,11pt]{amsart}

	\usepackage[utf8]{inputenc}
	\usepackage[T1]{fontenc}
	\usepackage[french]{babel}
	\usepackage[protrusion=true,expansion=true,kerning, spacing]{microtype}
	\microtypecontext{spacing=french}
	\usepackage{CJKutf8} 
	\usepackage{cmap} 

	\usepackage{xspace}
	
	\usepackage[textwidth=16cm,hcentering,heightrounded,foot=1cm]{geometry}
	\usepackage{caption}
 \usepackage{titlesec}
	\titleformat{name=\section}[hang]{\scshape \large \centering}{\thesection.~}{0em}{}[]

	\titleformat{\subsection}[hang]{\bfseries}{\thesubsection.~}{0em}{}[]

	\titleformat{\subsubsection}[runin]{\bfseries}{\thesubsubsection.}{.2em}{}[]

 \usepackage{titletoc}
\titlecontents{section}[2pc]{\addvspace{1pc}}
{\contentslabel[\textsc{\thecontentslabel}]{1pc}}
{}{\dotfill\contentspage}
[\addvspace{2pt}]

\titlecontents{subsection}[4pc]{}
{\contentslabel[\subsectionname\thecontentslabel]{2pc}}
{}{\hfill\contentspage}
[]

\usepackage{csquotes}
\usepackage[
bibencoding=utf8,%
giveninits=true,
bibstyle=ieee-alphabetic,
citestyle=alphabetic,%
sorting=nyt,%
language=french,%
isbn=false,%
doi=false,%
eprint=true,
backref=true, 
hyperref=true,
url=true,
minbibnames=4,
maxbibnames=4,
dashed=true,
backend=biber,
]{biblatex}
	\usepackage{amsthm,amssymb,amsmath}
	\usepackage{amssymb,amsfonts}
	\usepackage[fixamsmath]{mathtools} 
	\usepackage{accents} 
	\usepackage{tikz-cd}
		\numberwithin{equation}{subsection}

	\usepackage{graphics,graphicx}
	\usepackage{color}
	\usepackage{tikz}
	\usepackage{caption} 
	\usepackage{subcaption} 
	\usepackage{wrapfig}

\theoremstyle{definition}
	\newtheorem{defin}{Définition}[section]
	\newtheorem{ex}[defin]{Exemple}

\theoremstyle{definition}
	\newtheorem{rem}[defin]{Remarque}
	
\theoremstyle{plain}
	\newtheorem{theo}[defin]{Théorème}
	\newtheorem{prop}[defin]{Proposition}
	\newtheorem{cor}[defin]{Corollaire}
	\newtheorem{lem}[defin]{Lemme}

\newcommand{\C}{\mathbf{C}}
\newcommand{\A}{\mathbf{A}}

\newcommand{\N}{\mathbf{N}}
\newcommand{\R}{\mathbf{R}}
\newcommand{\Q}{\mathbf{Q}}
\newcommand{\Z}{\mathbf{Z}}

\newcommand{\F}{\mathbf{F}}

\newcommand{\G}{\mathbf{G}}

\newcommand{\OO} {{\mathcal O}} 
\newcommand{\mA} { \mathcal{A}}
\newcommand{\mF} {\mathcal{F}}
\newcommand{\mB} {\mathcal{B}}
\newcommand{\mD} {\mathcal{D}}

\newcommand{\mC} {\mathcal{C}}

\newcommand{\mL} {\mathcal{L}}
\newcommand{\mP} {\mathcal{P}}
\newcommand{\mfm} {\mathfrak{m}}
\newcommand{\mfM} {\mathfrak{M}}

\newcommand{\mfS} {\mathfrak{S}}

\newcommand{\Ker}{\operatorname{Ker}}

\newcommand{\End} {\operatorname{End} }
\newcommand{\Hom} {\operatorname{Hom} }

\newcommand{\Card} {\operatorname{Card} }
\newcommand{\Aut} {\operatorname{Aut} }
\newcommand{\Gal} {\operatorname{Gal} }

\newcommand{\GL}{\operatorname{GL}}
\newcommand{\oD}{\operatorname{D}}
\newcommand{\Gr}{\operatorname{Gr}}
\newcommand{\LGr}{\operatorname{LGr}}

\newcommand{\ppcm}{\operatorname{ppcm}}
\newcommand{\Rep}{\operatorname{Rep}}
\newcommand{\Fil}{\operatorname{Fil}}

\newcommand{\Sp}{\operatorname{Sp}}

\usepackage[pdftex,
hidelinks,
hyperindex=true, bookmarks=true, bookmarksopen=true, bookmarksnumbered=true,
pdfauthor={Séverin Philip},
pdftitle={Groupes de monodromie finie des variétés abéliennes},
pdfsubject={},
pdfstartview=FitH
]{hyperref}

\begin{document}
\author{Séverin Philip}
\address{Department of Mathematics, Stockholms universitet
SE-106 91 Stockholm, Sweden}
	\email{severin.philip@math.su.se}
 \thanks{L'auteur est soutenu par la bourse JSPS KAKENHI numéro 22F22015. L'auteur souhaite remercier C.~Cornut, A.~Tamagawa, G.~Rémond, A.~Mézard et B.~Collas pour les discussions fructueuses au sujet de ce travail.}

\title{Groupes de monodromie finie des variétés abéliennes \\[0.5em]
\textnormal{ Version du \today}}

\maketitle

\begin{abstract}
   Les groupes de monodromie finie des variétés abéliennes sont introduits par Grothendieck et représentent l'obstruction locale à la réduction semi-stable. Il est connu depuis Serre-Tate que ces groupes s'injectent dans les groupes d'automorphismes de variétés semi-abéliennes sur corps fini. On en donne dans ce texte une réciproque par une construction par déformation utilisant la théorie de Hodge $p$-adique entière. On obtient une caractérisation des groupes finis réalisables comme groupe de monodromie finie en dimension fixée ainsi qu'une application au degré de semi-stabilité. Précisément, on calcule une borne optimale, ne dépendant que de la dimension, au degré d'une extension finie fournie par le théorème de réduction semi-stable de Grothendieck.  
\end{abstract}

\tableofcontents

\section{Introduction}

\addtocounter{subsection}{1}

\subsubsection{} Pour une variété abélienne $A$ sur un corps de nombres ou $p$-adique $K$, Grothendieck montre dans \cite{sga} exposé IX que l'obstruction à la réduction semi-stable en une place non-archimédienne $v$ de $K$ est régie par un groupe fini. On note ce groupe $\Phi_{A,v}$ et on l'appelle groupe de monodromie finie de $A$ en $v$. On s'intéresse dans ce texte à la réalisabilité des groupes finis comme groupes de monodromie finie en dimension fixée et sur des corps de nombres. Cette question est étudiée par Silverberg et Zarhin dans \cite{SZ05} pour les surfaces abéliennes en caractéristique positive. Le résultat principal de ce texte donne une caractérisation de ces groupes par les automorphismes des variétés semi-abéliennes sur les corps finis. 

\begin{theo}\label{theo:main}
    Soit $G$ un groupe fini. Alors il existe une variété abélienne $A$ de dimension $g$ sur un corps de nombres $K$ et une place non-archimédienne $v$ de caractéristique résiduelle $p$ de $K$ telle que 
    \[
    \Phi_{A,v}= G
    \]
    si et seulement si, le groupe $G$ est le groupe de Galois d'une extension $L/K^{\mathrm{nr}}$ où $K^{\mathrm{nr}}$ est l'extension maximale non ramifiée d'un corps $p$-adique $K$, et s'il existe une injection $G\hookrightarrow \Aut(A_0, \lambda_0)$ pour une variété semi-abélienne polarisée $(A_0,\lambda_0)$ de dimension $g$ sur un corps fini de caractéristique $p$.  
\end{theo}


La notion de variété semi-abélienne polarisée est introduite en définition~\ref{def:vsapol}. Le sens direct du théorème est bien connu et utilisé dans différents articles tels que \cite{ST68} et \cite{SZ98}. La nouveauté est l'obtention de la réciproque, ce qui est fait en procédant par déformation grâce à la théorie de Hodge $p$-adique. 

\smallskip

Une application au degré de semi-stabilité en partie~\ref{sec:degsemistab} complète les résultats de \cite{Ph221} en remplaçant la construction par torsion galoisienne de variétés abéliennes avec multiplication complexe de \cite{Phi222} lorsque $p=2$. On calcule une borne optimale pour le degré minimal d'une extension sur laquelle une variété abélienne de dimension $g$ sur un corps de nombres à réduction semi-stable, donnée par la borne de Minkowski $M(2g)$, plus petit commun multiple des cardinaux des sous-groupes fini de $\GL_{2g}(\Q)$. Cela permet d'énoncer une version effective du théorème de réduction semi-stable de Grothendieck. 

\begin{theo}[de réduction semi-stable, forme effective] \label{theo:grotheff}
    Soient un entier $g\geq 1$ et une variété abélienne $A$ de dimension $g$ sur un corps de nombres $K$. Il existe une extension finie $L/K$ de degré au plus $M(2g)$ telle que $A_L$ a réduction semi-stable. De plus, il existe une variété abélienne $A$ de dimension $g$ sur un corps de nombres $K$ telle que pour toute extension finie $L/K$ de degré strictement plus petit que $M(2g)$ la variété abélienne $A_L$ n'a pas réduction semi-stable. 
\end{theo}

\subsubsection{} Avant de présenter la structure du texte et les différentes étapes des preuves on introduit les notions en jeu avec le sens direct du théorème~\ref{theo:main}. On considère donc une variété abélienne $A$ de dimension $g$ sur un corps de nombres $K$ et on fixe une place non-archimédienne $v$ de $K$ de caractéristique résiduelle $p$. 

\begin{defin}[\cite{sga}, exposé IX, partie 4] \label{def:monod}
    On note $\Phi_{A,v}$ et on appelle groupe de monodromie finie de $A$ en $v$ le groupe de Galois de la plus petite extension $L$ de $K_v^{\mathrm{nr}}$ sur laquelle $A$ a réduction semi-stable.
\end{defin}
Ce groupe est déterminé par la représentation $\ell$-adique de $A$ issue de son module de Tate, pour tout choix de $\ell\neq p$. En particulier, il est invariant par $K$-isogénie. De plus, la définition assure que $\Phi_{A,v}$ est un groupe de ramification en $p$, c'est-à-dire le groupe de Galois d'une extension totalement ramifiée de corps $p$-adiques comme dans l'énoncé~\ref{theo:main}. 

\medskip

 La variété abélienne $A_L$ est alors munie d'une donnée de descente canonique $(f_{\sigma})_{\sigma\in \Phi_{A,v}}$ qui commute à $\lambda_L$, où $\lambda$ est un choix de polarisation sur $A$. Par propriété de Néron cette donnée de descente, ainsi que $\lambda_L$, s'étendent aux modèles de Néron de $A_L$ et $A_L^{\vee}$ sur $\OO_L$. Comme l'extension $L/K_v^{\mathrm{nr}}$ est totalement ramifiée, le passage à la fibre spéciale du modèle de Néron de $A_L$ transforme la donnée de descente $(f_{\sigma})_{\sigma\in \Phi_{A,v}}$ en une action fidèle de $\Phi_{A,v}$ sur $A_0$, la composante neutre de la réduction de $A_L$. Cette action est alors compatible à la réduction $\lambda_0$ de $\lambda_L$. Le morphisme $\lambda_0$ est un morphisme de suites exactes
\[
\begin{tikzcd}
    0 \arrow[r] & T_0 \arrow[r] \arrow[d, "\lambda_{T_0}"] & A_0 \arrow[r, "p"] \arrow[d, "\lambda_0"] & B_0 \arrow[r]  \arrow[d, "\lambda_{B_0}"]& 0\\
    0 \arrow[r] & T_0' \arrow[r] & A_0' \arrow[r, "p"] & B_0^{\vee} \arrow[r] & 0
\end{tikzcd}
\]
où $A_0'$ est la réduction de $A_L^{\vee}$. Cette situation motive la définition suivante.

\begin{defin}\label{def:vsapol}
Un morphisme $\lambda_0$ de variétés semi-abéliennes $A_0\to A_0^t$ sur un corps fini $k$ 
\[
\begin{tikzcd}
    0 \arrow[r] & T_0 \arrow[r] \arrow[d, "\lambda_{T_0}"] & A_0 \arrow[r, "p"] \arrow[d, "\lambda_0"] & B_0 \arrow[r]  \arrow[d, "\lambda_{B_0}"]& 0\\
    0 \arrow[r] & T_0^t \arrow[r] & A_0^t \arrow[r, "p"] & B_0^{\vee} \arrow[r] & 0
\end{tikzcd}
\]
est une polarisation si le morphisme induit $\lambda_{T_0}$ est une isogénie et le morphisme $\lambda_{B_0}$ est une polarisation. Le degré $\deg \lambda_0$ de $\lambda_0$ est alors $m^2\cdot \deg \lambda_{B_0}$ où $m$ est l'ordre de $\Ker \lambda_{T_0}$. Une variété semi-abélienne sur un corps fini muni d'une polarisation est dite polarisée.
\end{defin}

La compatibilité de l'action de $\Phi_{A,v}$ avec la suite exacte de $A_0$ ainsi que $\lambda_0$ assure que l'on a une injection, à comparer avec le théorème~5.2 de \cite{SZ98},
\[
\iota\colon \Phi_{A,v} \hookrightarrow \Aut (T_0, \lambda_{T_0}) \times \Aut (B_0, \lambda_{B_0})
\]
où $\Aut (T_0, \lambda_{T_0})=\{ \alpha \in \Aut T_0 \mid \exists \alpha'\in \Aut T_0',~\lambda_{T_0} \alpha= \alpha' \lambda_{T_0}\}$ et $\Aut (B_0, \lambda_{B_0})$ est défini de manière classique $\Aut (B_0, \lambda_{B_0})= \{ \alpha\in \Aut B_0 \mid \lambda_{B_0} \alpha= {\alpha^{-1}}^{\vee} \lambda_{B_0} \}$ du fait que la monodromie finie préserve la polarisation $\lambda_{B_0}$ induite sur $B_0$. 

\begin{defin}
    On note $\Aut (A_0,\lambda_0)$ l'ensemble des automorphismes de $A_0$ dont l'image dans $\Aut T_0 \times \Aut B_0$ a pour première projection un élément de $\Aut (T_0, \lambda_{T_0})$ et pour deuxième projection un élément de $\Aut (B_0, \lambda_{B_0})$.
\end{defin}

On a finalement obtenu la condition nécessaire annoncée pour que $G$ soit réalisable comme groupe de monodromie finie. Le théorème~\ref{theo:main} assure que cette condition est, en fait, suffisante.

\begin{prop} \label{prop:sensdirect}
    Soit $G$ un groupe fini. Si $G$ est réalisable comme groupe de monodromie finie d'une variété abélienne de dimension $g$ sur un corps de nombres alors il existe une variété semi-abélienne polarisée $(A_0,\lambda_0)$ de dimension $g$ sur un corps fini et une injection $G\subset \Aut (A_0, \lambda_0)$.
\end{prop}

\subsubsection{} Le théorème~\ref{theo:resumconstruction}, obtenu en fin de partie~\ref{sec:constru}, donne la réciproque à la proposition~\ref{prop:sensdirect} dont notre résultat principal découle directement. Elle repose sur une caractérisation de la monodromie finie par descente galoisienne, objet du théorème~\ref{theo:descentemonod}. \'Etant donné une variété semi-abélienne polarisée $(A_0, \lambda_0)$ sur un corps fini munie d'une action fidèle d'un groupe $G$ de ramification, il suffit donc de trouver une variété abélienne munie d'une donnée de descente qui relève, en un sens à préciser, la donnée $(A_0,\lambda_0, G)$. En partie~\ref{sec:constru} on entreprend de faire cette construction par déformation en utilisant la théorie de Hodge $p$-adique entière, dont les différentes catégories semi-linéaires et foncteurs pertinents à notre propos sont rappelés en partie~\ref{sub:rappelsHodge}. La construction par déformation suit, en partie~\ref{sub:deformHodge}. Elle repose de façon essentielle sur une généralisation du théorème de déformation de Serre-Tate aux variétés semi-abéliennes par Bertapelle et Mazzari dans \cite{BM19} qui assure que la théorie de la déformation de telles variétés sur les corps finis est donnée par celle de leur groupe $p$-divisible. La théorie de Hodge $p$-adique entière permet alors de traiter des groupes $p$-divisibles sur l'anneau des entiers d'un corps $p$-adique. Plus précisément, d'après Breuil, Kisin et Kim dans \cite{Kis06} et \cite{Kim12} la catégorie $p-\mathrm{div}/\OO_K$ des groupes $p$-divisibles sur $\OO_K$ est équivalente à celle des modules de Breuil-Kisin pour un corps $p$-adique $K$ -- voir partie~\ref{sub:rappelsHodge} pour les définitions précises. On obtient, par l'existence de certaines filtrations admissibles sur les $\varphi$-modules issus de variétés semi-abéliennes, un schéma semi-abélien qui relève $A_0$, à interpréter comme l'extension de Raynaud de la variété abélienne finale. Une dégénérescence de celui-ci, suivant la théorie de Faltings et Chai aux chapitres 2 et 3 de \cite{FC90}, fournit la variété abélienne cherchée. 

\medskip

\subsubsection{} L'existence de filtrations admissibles sur les $\varphi$-modules, nécessaire à la construction présentée au paragraphe précédent, est le cœur technique de ce papier et occupe la partie~\ref{sec:phimod}. On introduit d'abord la notion de $\varphi$-module abélien et semi-abélien polarisé sur lesquelles on pose le problème de l'existence de filtration admissibles compatibles avec leur structure et vérifiant certaines propriétés. Précisément, les filtrations doivent être compatibles à la polarisation et produire une donnée de descente pour une extension galoisienne totalement ramifiée à partir d'une action fidèle de son groupe de Galois sur le $\varphi$-module -- voir le théorème~\ref{theo:EAdm} et le paragraphe le précédant. La preuve de l'existence de ces filtrations, grâce à une décomposition des $\varphi$-modules semi-abéliens munis d'action de groupes finis, est réduite à cette même question d'existence dans deux cas particuliers. D'abord, celui où le $\varphi$-module est supersingulier, c'est-à-dire isocline de pente $1/2$ et ensuite celui où il ne possède pas de composantes de pente $1/2$. Dans ce second cas, on montre l'existence d'un ouvert de Zariski dans l'espace des filtrations admissibles. Dans le premier, un tel ouvert n'existe pas toujours et on utilise l'action d'éléments appelés perturbateurs pour conclure. 

\subsubsection{} L'application principale du théorème~\ref{theo:main} concerne le degré de semi-stabilité des variétés abéliennes. Précisément, pour une variété abélienne $A$ sur un corps de nombres $K$ on note $d(A)$ le minimum des degrés des extensions finies $L/K$ telles que $A_L$ a réduction semi-stable. On considère alors l'entier $d_g$ défini comme minimum des entiers $d\geq 1$ vérifiant $d(A)\leq d$ pour toute variété abélienne $A$ de dimension $g$ sur un corps de nombres. Dans \cite{Ph221}, grâce à la méthode de torsion évoquée ci-dessus, les inégalités suivantes sont établies pour tout $g\geq 1$
\[
 \frac{M(2g)}{2^{g-1}} \leq d_g \leq M(2g)
\]
où $M(n)=\prod_{p\text{ premier}} p^{r(n,p)}$ avec $r(n,p)=\sum_{i\geq 0} \big\lfloor \frac{n}{p^i(p-1)} \big \rfloor$ est la borne de Minkowski. On montre, en partie~\ref{sub:maxd}, pour tout $g\geq 1$, l'égalité
\[
d_g=M(2g)
\]
en réalisant un $2$-Sylow du produit en couronne $Q_8\wr \mathfrak{S}_g$ comme groupe de monodromie finie en dimension $g$ grâce au théorème~\ref{theo:main} et cela donne donc le théorème~\ref{theo:grotheff}.

\medskip

En partie~\ref{sub:deploy} on étudie l'analogue $d_g^{d\acute{e}p}$ de $d_g$ où l'on demande de plus que les tores intervenant aux places de mauvaise réduction soient déployés. On montre ici, à nouveau pour $g\geq 1$, l'égalité plus forte
\[
d_g^{d\acute{e}p}=M(2g).
\]

\section{Isocristaux de variétés semi-abéliennes et filtrations admissibles} \label{sec:phimod}

\subsection{Isocristaux, filtrations admissibles et grassmanniennes} \label{sub:isoc}

\subsubsection{} Soient $\overline{k}$ la clôture algébrique d'un corps fini de caractéristique $p$ et $K_0$ le corps des fractions de l'anneau des vecteurs de Witt $W(\overline{k})$, c'est-à-dire $K_0\simeq \Q_p^{\mathrm{nr}}$. Un $\varphi$-module, ou isocristal, $D$ sur $K_0$ est un $K_0$-espace vectoriel $D$ de dimension finie muni d'un opérateur bijectif $\varphi$ semi-linéaire. Précisément, si l'on note $\sigma$ le Frobenius sur $W(k),$ la semi-linéarité s'exprime alors, pour tout $\lambda\in K_0$ et $x\in D$, par l'égalité $\varphi(\lambda x)= \sigma(\lambda)\varphi(x)$. On note $\operatorname{MF}^{\varphi}$ la catégorie des $\varphi$-modules, c'est une sous-catégorie pleine de celle des $\varphi,N$-modules $\operatorname{MF}^{\varphi,N}$ décrite en définition~8.13 de \cite{FO22}. 

Plusieurs invariants des $\varphi$-modules vont intervenir dans cette partie, dont en particulier leur pente. Le premier, noté $t_N(D)$ pour le $\varphi$-module $D$, peut se définir à partir d'une matrice $A$ représentant l'opérateur semi-linéaire $\varphi$ -- voir les définitions 8.18, 8.20 de \cite{FO22} -- en posant $t_N(D)=v_p(\det A)$. La pente se définit alors comme le quotient $t_N(D)/\dim_{K_0} D$. Par le théorème de Dieudonné-Manin -- voir le théorème~8.25 de \cite{FO22} -- le $\varphi$-module $D$ admet une décomposition isocline
\[
D= \bigoplus_{\mu\in \Q} D_{\mu}
\]
où $D_{\mu}$ est un $\varphi$-module pur de pente $\mu\in \Q$ -- voir la définition~8.23 de \cite{FO22}. L'intérêt principal de cette notion dans ce texte est que tout sous-$\varphi$-module d'un module pur de pente $\mu$ est encore pur de pente $\mu$, ce qui se déduit aisément de la décomposition semi-simple des $\varphi$-modules sur $K_0$.

\medskip

Pour une extension finie $K/K_0$, on considère dans ce texte des filtrations sur l'espace vectoriel $D_K$ à un seul cran : formellement de la forme $D_K\supset F \supset \{0\}$ -- $\Fil^0=D_K$, $\Fil^1=F$ et $\Fil^2=\{0\}$ dans les notations de \cite{FO22}. Une telle filtration sera confondue avec le sous-espace vectoriel $F\subset D_K$. La notion de filtration admissible pour $D$, donnée en partie~8.2.5 de \cite{FO22}, est ici celle d'un sous-espace vectoriel $F\subset D_K$ qui vérifie, pour tout sous-$\varphi$-module $N\subset D$,
\[
\dim N_K\cap F\leq \sum\limits_{\mu\in \Q} \mu d_{\mu}
\]
où $d_{\mu}$ est la dimension de la composante de pente $\mu$ de $N$ et avec égalité pour $N=D$. En effet, la notion d'admissibilité consiste à comparer l'invariant $t_N$ introduit précédemment à un autre $t_H$ -- voir la définition~8.36 et la proposition~8.37 de \cite{FO22} -- sur $D$ et ses sous-$\varphi$-modules munis de la filtration induite par l'intersection. Précisément -- voir la définition 8.39 de \cite{FO22} -- on demande l'inégalité $t_H(N)\leq t_N(N)$ pour tout sous-$\varphi$-module $N$ de $D$ avec égalité lorsque $N=D$. Par la décomposition isocline on a l'égalité $t_N(N)=\sum_{\mu} \mu d_{\mu}$ et comme nos filtrations ont un seul cran il est immédiat que $t_H(N)=\dim N_K\cap F$.

\smallskip

On note $\operatorname{MF}^{\varphi,0,1}_K$ la catégorie des $\varphi$-modules filtrés admissibles sur $K$ dont les filtrations ont un seul cran, c'est une sous-catégorie pleine de la catégorie $\operatorname{MF}^{\varphi,N}_K$ des $\varphi,N$-modules filtrés admissibles introduite en définition~8.39 de \cite{FO22}. On note les objets de la catégorie $\operatorname{MF}^{\varphi,0,1}_K$ sous la forme d'un couple $(D,F)$ où $D$ est un $\varphi$-module et $F\subset D_K$ une filtration admissible. Les morphismes dans cette catégorie sont les applications linéaires $f\colon (D_1,F_1)\to (D_2,F_2)$ compatibles aux Frobenius $\varphi$ et vérifiant $f(F_1)\subset F_2$. On dispose d'une dualité sur cette catégorie. Le dual d'un $\varphi$-module filtré $(D,F)$ est donné par $(D^{\vee}, F^{\perp})$ où $D^{\vee}$ est muni d'une structure de $\varphi$-module par 
\[
\varphi(f)=\sigma \circ f \circ \varphi^{-1}
\]
où $f\in D^{\vee}$. Dans le cas où $(D,F)$ est autodual il est de dimension paire $2g$. En effet, les pentes de $D$ apparaissent alors par paires $(\mu, 1-\mu)$, le dual d'un $\phi$-module $D_{\mu}$ pure de pente $\mu$ étant $D_{1-\mu}$ et si $\mu=1/2$ alors $D_{\mu}$ est de dimension paire. Il suit de plus que la filtration $F$ est un sous-espace de dimension $g$ de $D_K$. Dans cette situation, le foncteur $K\mapsto \{ F\subset D_K \mid F \text{ est une filtration admissible}\}$ est représentable par un ouvert analytique $\mF^{\mathrm{adm}}$ au sens de Berkovic de la grassmannienne $\mathrm{Gr}_{2g,g}$ des sous-espaces de dimension $g$ de $D$ de complémentaire un compact, voir \cite{DOR10} Proposition~8.2.1. Cette description, qui est valable en plus grande généralité, n'est pas utile dans ce texte mais la notation $\mF^{\mathrm{adm}}$ est employée.


\subsubsection{}\label{subsub:galoisextexist} On introduit maintenant la situation précise qui va nous intéresser pour le reste de cette partie. 

\begin{defin}
    Un $\varphi$-module $D$ est dit abélien s'il est autodual et si toutes ses pentes sont dans $\Q\cap [0,1]$. On dira de plus qu'il est polarisé s'il est muni d'un isomorphisme $\lambda_0$ sur son dual qui détermine une forme bilinéaire alternée sur $D$. 

    Un $\varphi$-module $D$ est dit semi-abélien s'il s'inscrit dans une suite exacte
    \[
    \begin{tikzcd}
        0 \arrow[r] & D_T \arrow[r] & D \arrow[r] &  D_B \arrow[r] & 0
    \end{tikzcd}
    \]
    où $D_T$ est pur de pente $1$ et $D_B$ est un $\varphi$-module abélien. On dira de plus qu'il est polarisé s'il existe un $\varphi$-module semi-abélien $D^t$ et un isomorphisme $\lambda_0\colon D\to D^t$ qui respecte les suites exactes définissant $D$ et $D^t$ et qui munit $D_B$ d'une structure de $\varphi$-module abélien polarisé. En particulier, $D^t_B$ est le dual de $D_B$.
\end{defin}


Soit $G$ un groupe fini de ramification. Par définition, il existe une extension galoisienne $L/K$ de groupe de Galois le groupe opposé $G^{\mathrm{op}}$ de $G$ où $K/K_0$ est modérément ramifiée. Le groupe $G^{\mathrm{op}}$ agit donc naturellement à gauche sur $L$ et ainsi $G$ agit naturellement à droite sur $L$. On considère dans la suite cette action naturelle à droite de $G$ et on l'appelle l'action galoisienne de $G$. 

\bigskip


Soit alors un $\varphi$-module semi-abélien polarisé $D$ muni d'une action fidèle du groupe fini $G$, l'action de $G$ est toujours considérée compatible avec la structure de $\varphi$-module semi-abélien polarisé. Le groupe $G$ agit alors naturellement de deux manières sur l'espace vectoriel $D_L$. La première action est linéaire et provient de l'inclusion $G\subset \Aut_{\varphi} D$. La deuxième se déduit de l'action galoisienne sur $L$ et se fait sur le deuxième facteur de $D_L=D\otimes_{K_0} L$. On notera la première action $\cdot_{\mathrm{lin}}$ et la deuxième $\cdot_{\mathrm{gal}}$. Ces deux actions l'une à gauche et l'autre à droite commutent entre elles et on pourra considérer l'action diagonale donnée par $h\otimes h^{-1}$ pour $h\in G$. On va chercher des filtrations qui vérifient une propriété de compatibilité à cette action diagonale. Précisément, une condition de stabilité de $F$ pour l'action diagonale de $G$ introduite précédemment, ou encore que les actions galoisiennes et linéaires de $G$ sur $F$ coïncident, ce qui s'écrit pour tout $h\in G$
\[
h\cdot_{\operatorname{gal}} F= h\cdot_{\operatorname{lin}} F.
\]

\subsubsection{} Le reste de la partie~\ref{sec:phimod} est dédié à l'étude de l'existence de filtrations admissibles dont les propriétés sont données au paragraphe précédent et résumées dans l'énoncé du théorème~\ref{theo:EAdm}. Les sous-parties suivantes établissent les résultats préliminaires à la démonstration du théorème~\ref{theo:EAdm} et donnent différents critères pour l'existence de filtrations admissibles avec les propriétés voulues. 

\begin{theo} \label{theo:EAdm}
     Soit $D$ un $\varphi$-module semi-abélien polarisé sur $K_0$ muni d'une action d'un groupe fini $G$ de ramification. Alors, pour toute extension galoisienne $L/K$ de groupe $G$ avec $K/K_0$ modérément ramifiée et finie, il existe une filtration $F\subset D_L$ qui vérifie les propriétés suivantes : 
\begin{enumerate}
    \item \label{enum:lag} $F/(F\cap (D_{T_0})_L)$ est lagrangien pour $\lambda_{B_0}$,
    \item \label{enum:adm}$F$, $F\cap (D_{T_0})_L$ et $F/(F\cap D_{T_0})$ sont des filtrations admissibles,
    \item \label{enum:stabgal} $F$ est stable sous l'action diagonale de $G$, autrement dit l'orbite de $F$ sous l'action linéaire de $G$ correspond naturellement avec son orbite sous l'action galoisienne, ce qui s'écrit
    
\[
\forall h\in G,~ h\cdot_{\mathrm{lin}} F= h \cdot_{\mathrm{gal}} F.
\]
\end{enumerate}

\end{theo}

Rappelons que, sauf mention du contraire, l'on s'est fixé une telle extension arbitraire $L/K$ au paragraphe~\ref{subsub:galoisextexist} et que toute extension $L/K$ apparaissant dans la suite seront prises au-dessus de $K_0$. 

\smallskip

Montrons tout d'abord que pour prouver le théorème il suffit de le faire pour dans le cas où $D$ est abélien.

\begin{prop}\label{prop:reductab}
    Soit $D$ un $\varphi$-module semi-abélien polarisé muni d'une action d'un groupe fini $G$ de ramification. Alors la conclusion du théorème~\ref{theo:EAdm} vaut pour $D$ si elle vaut pour $D_{B}$ muni de l'action de $G'$, quotient de $G$, déduite de la projection $\pi\colon D\to D_B$.
\end{prop}
\begin{proof}
    Soit $L/K$ une extension galoisienne de groupe $G$ dont la sous-extension $L'/K$ est galoisienne de groupe $G'$. Soit $F'$ une filtration de $(D_{B_0})_{L'}=(D/D_{T_0})_{L'}$ avec les propriétés (\ref{enum:lag}), (\ref{enum:adm}) et (\ref{enum:stabgal}) pour l'action de $G'$. Il suffit de vérifier que le relevé $F= \pi^{-1}(F')$ convient. La propriété (\ref{enum:lag}) pour $F$ portant sur $F'$, elle est vérifiée par hypothèse. Il en va de même pour la propriété d'admissibilité dans le quotient. La propriété (\ref{enum:adm}) demande l'admissibilité de $F\cap D_{T_0}$ ce qui se traduit par l'inclusion $D_{T_0} \subset F$ et du fait que $D_{T_0}$ est pur de pente $1$ la question de l'admissibilité de $F$ pour $D$ se fait par un calcul direct. Il ne reste donc qu'à vérifier la propriété (\ref{enum:stabgal}). Cela se fait sans encombre du fait que le morphisme $\pi$ commute aux actions linéaires et galoisiennes de $G$. Soit $h\in G$, d'image $h'$ dans $G'$, l'égalité $h'\cdot_{\mathrm{lin}} F'=h'\cdot_{\mathrm{gal}} F'$ assure que
    \begin{align*}
        h\cdot_{\mathrm{lin}}F &= h\cdot_{\mathrm{lin}}(\pi^{-1}(F')) \\
        &= \pi^{-1}(h'\cdot_{\mathrm{lin}} F') \\
        &= \pi^{-1}(h'\cdot_{\mathrm{gal}} F')\\
        &= h\cdot_{\mathrm{gal}} (\pi^{-1}(F')).
    \end{align*}
\end{proof}

 Suivant ce résultat on va se restreindre dans toute la suite de cette partie au cas où $D$ est abélien muni d'une polarisation $\lambda_0$. La grassmannienne lagrangienne $\operatorname{LGr}_{2g,g}$ pour $\lambda_0$ est alors le fermé de $\mathrm{Gr}_{2g,g}$ dont les $K$-points sont les sous-espaces lagrangiens de $(D_K,\lambda_0)$ pour un corps $K/K_0$. 

\subsubsection{} On termine cette sous-partie par l'étude des filtrations qui vérifient les propriétés (\ref{enum:lag}) et (\ref{enum:stabgal}) dans la situation introduite au paragraphe~\ref{subsub:galoisextexist} en supposant de plus $D$ abélien. On aura besoin pour cela d'un lemme classique d'existence de bases invariantes par une action galoisienne semi-linéaire sur des espaces vectoriels, à comparer avec le lemme~II.5.8.1 de \cite{Si09}.
\begin{lem}\label{lem:vectgalois}
    Soient $L/K$ une extension finie galoisienne de groupe $G$ et $E$ un $L$-espace vectoriel muni d'une action semi-linéaire du groupe fini $G$. Alors l'espace $E^G$ des éléments de $E$ fixés par l'action semi-linéaire de $G$ sur $E$ est un $K$-espace vectoriel et il vérifie $E_L^G\otimes_K L=E_L$. En particulier, c'est un $K$-espace vectoriel de dimension $\dim_L E$.
\end{lem}
\begin{proof}
    Soient $(\alpha_i)_{i\in \{1,\dots ,n\}}$ une base de $L/K$ et $v\in E$. On montre que $v$ s'obtient comme combinaisons linéaires à coefficients dans $L$ de vecteurs de $E^G$. On note $v_i$, pour $1\leq i \leq n$ le vecteur ligne de la matrice produit $[g(\alpha_i)]_{i,g\in G\times \{1,\dots, n\}}{[g\cdot v]^{\mathrm{T}}}_{g\in G}$. Par définition le vecteur 
    \[
    v_i=\sum\limits_{g\in G} g(\alpha_i)g\cdot v
    \]
    est invariant par $g$ pour tout $i$. La famille $(v_i)_{i\in \{1,\dots, n\}}$ est dans $E_L^G$ et, du fait que la matrice $[g(\alpha_i)]_{i,g\in G\times \{1,\dots, n\}}$ est inversible, elle engendre les vecteurs $g\cdot v$ sur $L$. 
\end{proof}

\begin{prop}\label{prop:descpointsgrassmann}
    Une filtration $F\subset D_L$ vérifie la propriété~(\ref{enum:stabgal}) de stabilité galoisienne si et seulement si $F=F^G\otimes_K L$ avec $F^G\subset D^G_L$ sous-$K$-espace vectoriel de dimension $g$. En particulier, les telles filtrations de $D_L$ sont données par les points $K$-rationnels de la grassmannienne $\Gr_{g}(D_L^G)$.  

    De plus, la correspondance analogue pour les filtrations lagrangiennes pour $\lambda_0$ a lieu. 
\end{prop}
\begin{proof}
    D'après le lemme~\ref{lem:vectgalois}, une filtration $F\subset D_L$ qui est stabilisés par l'action diagonale de $G$ admet une base de vecteurs invariants par $G$ ce qui donne l'équivalence annoncée. 

    \smallskip

    Toujours d'après le lemme~\ref{lem:vectgalois} $D_L^G$ est un $K$-espace vectoriel de dimension $2g$. Les sous-espaces vectoriels de dimension $g$ de ce dernier sont donc donnés par les $K$-points de la grassmannienne $\operatorname{Gr}_{2g,g}$. 

    \smallskip

    Pour en déduire l'assertion sur les filtrations lagrangiennes il suffit de vérifier que ${\lambda_0}_{|D_L^G}$ définit une forme bilinéaire alternée non dégénérée et que $\lambda_0={\lambda_0}_{|D_L^G}\otimes_{K} L$. La première assertion vient de ce que $\lambda_0$ est défini sur $K_0$ et que l'action linéaire de $G$ se fait par automorphismes symplectiques pour $\lambda_0$. La deuxième en suit directement.
\end{proof}

On va donner ici un critère pour qu'un sous-espace de la grassmannienne lagrangienne $\LGr$ contienne des filtrations vérifiant la propriété~(\ref{enum:stabgal}). Celui-ci sera utile pour certains cas particuliers où l'on vérifie la validité du théorème~\ref{theo:EAdm} dans les parties suivantes.  

\begin{prop} \label{prop:condit2Zariski}
    Tout ouvert de Zariski non vide de $\operatorname{LGr}_{2g,g}$ contient une filtration vérifiant la propriété (\ref{enum:stabgal}) de stabilité par l'action diagonale. 
\end{prop}
\begin{proof}
    L'isomorphisme $D_L^G\otimes_{K} L \simeq D_L$ donne un isomorphisme de variétés algébriques $\alpha\colon {\Gr_{2g,g}}_L\simeq \Gr_g (D_L^G)_L$. Soit alors $U\subset \Gr_{2g,g}$ un ouvert de Zariski. D'après la proposition~\ref{prop:descpointsgrassmann} un point $x\in U(L)$ vérifie la propriété (\ref{enum:stabgal}) si et seulement s'il est l'image d'un point $K$-rationnel de $\Gr_g(D_L^G)$ par l'isomorphisme $\alpha$. Comme la grassmannienne lagrangienne admet un recouvrement affine par des ouverts isomorphes à $\A^{\frac{g(g+1)}{2}}_{K}$, les points $K$-rationnels sont denses. En particulier un tel point existe dans $U$. 
\end{proof}

On déduit de  ces considérations un premier critère, qui ne sera pas toujours satisfait, d'existence de filtrations vérifiant les trois propriétés. 

\begin{cor}\label{cor:conditzariskEAdm}
S'il existe un ouvert de Zariski non vide $U$ de $\LGr_{\lambda_0}$ tel que $U(L)\subset \mF^{\mathrm{adm}}(L)$ pour une extension $L/K$ galoisienne de groupe $G$ alors il existe une filtration dans $U$ vérifiant les propriétés (\ref{enum:lag}), (\ref{enum:adm}) et (\ref{enum:stabgal}) avec l'extension $L/K$.
\end{cor}

\subsection{Cas supersingulier} \label{sub:supersing}

\subsubsection{} \label{subsub:ellipticsupersing}
On reprend la situation décrite au paragraphe~\ref{subsub:galoisextexist} et on suppose de plus que $D$ est abélien supersingulier au sens suivant.
\begin{defin}
    On dit qu'un $\varphi$-module $D$ abélien est supersingulier s'il est isocline de pente $1/2$. 
\end{defin}
Il suit de la définition que tous les sous-$\varphi$-modules de $D$ sont encore supersinguliers, donc de dimension paire et isocline de pente $1/2$. On note $2g$ la dimension de $D$.

\medskip

Le but de cette partie est de donner un critère pour l'existence de filtrations admissibles avec les propriétés (\ref{enum:lag}), (\ref{enum:adm}) et (\ref{enum:stabgal}) dans le cas supersingulier.

\medskip

La condition d'admissibilité sur $D$ s'exprime ici de façon simple du fait que $D$ est isocline de pente  $1/2$. Pour un sous-$\varphi$-module $N\subset D$ la condition s'écrit $\dim F\cap N_L \leq 1/2\dim N$. On montre dans un premier temps que l'espace analytique $\mF^{\mathrm{adm}}\subset \LGr_{2g,g}$ contient beaucoup de points fermés dans ce cas.

\begin{lem}\label{lem:supersingcondouvert} Soient $K'$ une extension finie non ramifiée de $\Q_p$ telle qu'il existe un $\varphi$-module $D'$ sur $K'$ avec $D=D'\otimes_{K'} K_0$ et $L'/K'$ une extension galoisienne totalement ramifiée telle que $L'/\Q_p$ est galoisienne. Soient $\tau\in \Gal (L'/\Q_p)$ un relèvement du générateur $\sigma$ de $\Gal(K'/\Q_p)$ à $L'$ et $\varphi_{\tau}\colon D'_{L'}\rightarrow D_{L'}$ l'extension $\tau$-semi-linéaire de $\varphi$ à $D'_{L'}$. Un sous-espace $F\subset {D'}_{L'}$ de dimension $g$ tel que $F\cap \varphi_{\tau}(F)=\{0\}$ définit une filtration admissible. De plus, l'ensemble des tels sous-espaces est donné par les $L'$-points d'un ouvert de Zariski non vide de la restriction de Weil de la grassmannienne $\LGr_{2g,g}$.
\end{lem}
\begin{proof}
    Soit $N\subset D'$ un sous-$\varphi$-module. Comme $N_{L'}$ est stable par $\varphi_{\tau}$ qui est inversible, on a $\varphi_{\tau}(F)\cap N_{L'}= \varphi_{\tau}(F)\cap \varphi_{\tau}(N_{L'})= \varphi_{\tau}(F\cap N_{L'})\subset N_{L'}$. Il suit $\dim \varphi_{\tau}(F\cap N_{L'}) = \dim F\cap N_{L'}$ et on obtient 
    \[
    \dim F\cap N_{L'}\leq \frac{1}{2}\dim N_{L'}
    \]
    du fait que $F\cap N_{L'}$ et $\varphi_{\tau}(F)\cap N_{L'}$ sont en somme directe et de même dimension dans $N_{L'}$. 

\medskip
    
    Pour la deuxième assertion, on remarque d'un côté que la compatibilité de $\varphi$ avec $\lambda_0$ assure que $\varphi_{\tau}(F)$ est lagrangien et d'un autre que la condition $F\cap \varphi_{\tau}(F)=\{0\}$ s'écrit sur les ouverts affines de $\operatorname{Gr}_{2g,g}/K'$ comme la non-annulation d'un déterminant non nul et, par restriction de Weil, les équations obtenues sont algébriques. 
\end{proof}

\begin{cor}\label{cor:supersingdense} Soit $L/K_0$ une extension finie telle que $L/\Q_p$ est galoisienne. L'ouvert $p$-adique des $L$-filtrations admissibles de l'espace topologique $\Gr_{2g,g}(L)$ est dense.
\end{cor}
\begin{proof}
    On note à nouveau $K'$ une extension finie non ramifiée de $\Q_p$ telle que $L$ descend en une extension galoisienne totalement ramifiée de $L'$ de $K'$ et telle que $D$ s'obtient comme changent de base d'un $\varphi$-module $D'$ sur $K'$.

    La restriction de Weil induit un isomorphisme entre les espaces topologiques $\Gr_{2g,g}(L')$ et $\operatorname{Res}_{L'/\Q_p}(\Gr_{2g,g})(\Q_p)$. D'après le lemme~\ref{lem:supersingcondouvert} il existe un ouvert de Zariski $U\subset \operatorname{Res}_{L'/\Q_p}(\Gr_{2g,g})$ dont les $\Q_p$-points, qui forment un ouvert dense, correspondent à des filtrations admissibles. Le résultat s'en déduit alors du fait que $K_0$, en tant que corps topologique est limite inductive de tels extensions $K'/\Q_p$. 
\end{proof}

Ce résultat, bien qu'utile, ne permet toutefois pas de montrer l'existence de filtrations vérifiant les trois propriétés voulues en général. 

\subsubsection{} Pour traiter pleinement de la situation introduite au paragraphe~\ref{subsub:galoisextexist} et dans le cas supersingulier il convient d'utiliser l'action de $G$.

\begin{theo} \label{theo:EAdmss}
Soit $D$ un $\varphi$-module supersingulier polarisé sur $K_0$. Soit $G$ un groupe de ramification fini qui agit linéairement sur le $K_0$-espace vectoriel $D$. S'il existe un élément $h\in G$ tel que la condition $\dim h(F)\cap F \leq 1$ définit un ouvert de Zariski non vide de la grassmannienne lagrangienne alors il existe une filtration sur $D$ vérifiant les propriétés (\ref{enum:lag}), (\ref{enum:adm}) et (\ref{enum:stabgal}).
\end{theo}

\begin{proof}
La condition d'admissibilité pour une filtration $F\subset D_L$ s'écrit comme avant $\dim F\cap N_L\leq \frac{1}{2} \dim N_L$ pour tout sous-$\varphi$-module $N\subset D$. 

\medskip

Soient $L/K$ une extension galoisienne de groupe $G$ et $h\in G$ tel que la condition $\dim F\cap h(F)\leq 1$ détermine un ouvert de Zariski non vide $U$ de la grassmannienne lagrangienne $\LGr_{\lambda_0}$ (un tel $h$ existe par hypothèse). Par la proposition~\ref{prop:condit2Zariski} il existe une filtration $F\in U(L)$ telle que les orbites induites par les actions galoisiennes et linéaires de $G$ sur $F$ coïncident. On a en particulier l'égalité
\[
h\cdot_{\mathrm{lin}} F= h\cdot_{\mathrm{gal}} F.
\]

Montrons que $F$ est admissible. Soit $N\subset D$ un sous-$\varphi$-module. Comme $N$ est défini sur $K_0$ il est stable par l'action galoisienne de $h$. On en déduit $h\cdot_{\mathrm{gal}}(F\cap N_L) \subset N_L$ et par hypothèse $\dim h\cdot_{\mathrm{gal}} (F\cap N_L) \cap (F\cap N_L)\leq 1$. Comme on a de plus $\dim h\cdot_{\mathrm{gal}} (F\cap N_L)= \dim (F\cap N_L)$, on obtient 
\[
2 \dim F\cap N_L - \dim  h\cdot_{\mathrm{gal}} (F\cap N_L) \cap (F\cap N_L) \leq \dim N
\]
et par parité de la dimension de $N$, $\dim F\cap N_L \leq \frac{1}{2} \dim N $.
    
\end{proof}

\begin{rem} \phantom{v}
\begin{itemize}
    \item[$(i)$] Il faut remarquer que la condition $\dim h(F)\cap F \leq 1$ définit toujours un ouvert de Zariski mais qu'il est vide dans de nombreux cas. 
    
    \item[$(ii)$] L'existence de $F$ vérifiant les trois propriétés découle sans avoir besoin que l'action de $G$ soit compatible à la structure de $\varphi$-module. Il n'est alors pas assuré qu'il existe des morphismes entre les $\varphi$-modules filtrés $(D,\sigma(F))$ pour $\sigma\in \Gal (L/K)$.
\end{itemize}   
\end{rem}

\begin{ex} \label{ex:Q8Sg}
    On donne ici une famille d'exemples où la condition du théorème~\ref{theo:EAdmss} est vérifiée. On considère le $\varphi$-module supersingulier $D^g_{1/2}$ pour $g\geq 1$ et où $D_{1/2}$ est le $\varphi$-module simple sur $\Q_2$ de dimension $2$ défini par 
    \[
    \varphi= \begin{pmatrix}
        0 & 2 \\
        1 & 0
    \end{pmatrix}.
    \]

    L'algèbre des endomorphismes $\End D_{1/2}$ est l'algèbre des quaternions classique $\mathbf{H}$ de base $(1,i,j,k)$ avec $i^2=j^2=-1$ et $ij=-ji=k$. Il suit $\End D_{1/2}^g= \mathrm{M}_{g}(\mathbf{H})$. On considère sur $D_{1/2}^g$ la forme bilinéaire alternée standard $\lambda_0^g$ où
    \[
    \lambda_0=\begin{pmatrix}
        0 & -1 \\  
        1 & 0 
    \end{pmatrix}
    \]
    qui est compatible à $\varphi$. 

    \smallskip

    Le produit en couronne $Q_8\wr \mfS_g$ se réalise matriciellement dans l'algèbre $\mathrm{M}_{g}(\mathbf{H})$ par les matrices de permutation à coefficients dans le groupe $Q_8$ qui sont, en particulier, orthogonales pour $\lambda_0$. La $2$-partie du cardinal de ce groupe est donnée par $r(2g,2)$ où $r(n,p)$ est la fonction rappelée en introduction donnant la puissance de $p$ dans la valeur de $M(n)$.

    \medskip
    
    On note $G_g$ un $2$-Sylow de ce groupe, qui par le lemme~3.5 de \cite{Phi222} est un groupe de ramification, et on vérifie l'existence d'un élément $h\in G$ vérifiant la condition du théorème~\ref{theo:EAdmss}. Pour ce dernier point, on considère l'élément $h$ donné par la matrice diagonale ayant l'élément $k\in Q_8$ sur ses coefficients diagonaux. Pour $g=1$ cela revient à prendre $h=k\in Q_8$. La condition du théorème~\ref{theo:EAdmss} se vérifie en considérant $h$ comme un élément de $\operatorname{M}_{2g}(\Q_4)$ avec son action linéaire sur $\Q^{2g}_4$. On peut passer à $\overline{\Q_4}$ pour vérifier que la condition donne un ouvert non vide pour la topologie de Zariski. Pour $g=1$, la matrice $2\times 2$ de $h$ est diagonalisable avec $2$ valeurs propres distinctes, notons-les $\alpha$ et $\beta$. Sur un ouvert affine de la grassmannienne, la condition s'écrit alors
    \[
    \begin{vmatrix}
        x & \alpha x\\
        y & \beta y
    \end{vmatrix} \neq 0
    \]
    qui donne bien un ouvert non vide $U$. Pour passer à la grassmannienne lagrangienne on utilise le lemme~\ref{lem:critlag} ci-après en remarquant que $U$ s'obtient comme l'intersection des ouverts déterminés par $F\cap M=\{0\}$ pour $M$ chacune des droites propres de $h$. On vérifie de la même façon que l'ouvert correspondant est non vide pour $g\geq 1$ en utilisant le fait que la polarisation considérée est la polarisation produit. 
\end{ex}

On termine par un lemme qui donne un critère pour vérifier l'hypothèse du théorème~\ref{theo:EAdmss}. On donne un autre critère en fin de partie~\ref{sub:kelem}.

\begin{lem}\label{lem:critlag}
    Soit $(V,\lambda)$ un espace symplectique. Alors si $M\subset V$ est un sous-espace de dimension plus petite que $1/2 \dim V$ alors il existe un sous-espace $F\subset V$ lagrangien pour $\lambda$ et vérifiant $F\cap M=\{0\}$. 
\end{lem}
\begin{proof}
    On montre le résultat par récurrence sur la dimension de $M$. Si $M$ est de dimension $0$ ou $1$ alors l'existence de lagrangiens transverses conclut. Sinon, soit $x\in M\setminus\{0\}$. Comme $\dim V\geq 2\dim M\geq \dim M +1$, il existe $y\in V \setminus \big(M\cup \langle x \rangle^{\perp}\big)$. On considère alors $V'= \langle x,y \rangle^{\perp}$ et $M'=V'\cap M\oplus \langle y\rangle$. Alors $V'$ est un espace symplectique qui vérifie $\dim V'= \dim V-2$ et $M'\subset V'$ vérifie l'hypothèse du fait que $M'\oplus \langle x,y \rangle= M\oplus \langle y \rangle $ d'où $\dim M'\leq \dim M-1$. L'hypothèse de récurrence appliquée au couple $(V',M')$ fournie un lagrangien $F'\subset V'$ tel que $F'\cap M'=\{0\}$. Alors $F=F'\oplus \langle y \rangle$ convient. L'égalité $F\cap M=\{0\}$ découle de $F'\cap (M\oplus \langle y \rangle= \{0\}$ et 
    \[
    F^{\perp}= (F')^{\perp}\cap \langle y \rangle^{\perp} = F'\oplus \langle x,y \rangle \cap V'\oplus \langle y \rangle = F'\oplus \langle y \rangle = F.
    \]
\end{proof}

Bien que le lemme précédent soit utile dans différentes situations, comme l'exemple~\ref{ex:Q8Sg} ou le cas de pentes différentes de $1/2$ de la partie suivante, on a besoin d'un autre critère pour vérifier la condition du théorème~\ref{theo:EAdmss} en général. Celui-ci utilise la notion d'endomorphisme perturbateur introduite en partie~\ref{sub:kelem} et est démontré en fin de cette partie.

\subsection{Isocristaux de pentes différentes de $1/2$}

\subsubsection{} On reprend à nouveau la situation du paragraphe~\ref{subsub:galoisextexist} en supposant ici que $D$ est abélien et n'admet pas de composante isocline de pente $1/2$. On va étudier les lagrangiens issues de la décomposition en composantes isoclines dans le but de donner un critère d'existence de filtrations avec les propriétés (\ref{enum:lag}), (\ref{enum:adm}) et (\ref{enum:stabgal}). Comme précédemment, on a l'écriture
\[
D=\bigoplus\limits_{\mu \in \mP} D_{\mu}
\]
où $\mP$ est l'ensemble des pentes de $D$. Par hypothèse, $1/2$ n'est pas une pente de $D$ et on peut donc coupler les pentes de $\mP$ par paire $(\mu, 1-\mu)$. On considère une base $(e_1,\dots, e_n)$ adapté à cette décomposition de $D$ et $(e_1^*, \dots, e_n^*)$ la base duale associée sur $D^{\vee}\simeq D$. Par la dualité pour les $\varphi$-modules, si $(e_j,\dots, e_k)$ est une base de la composante $D_{\mu}$ alors $(e_j^*, \dots , e_k^*)$ est une base de $D_{1-\mu}$ comme sous-espace de $D^{\vee}$. Par ailleurs, l'isomorphisme $\lambda_0\colon D\rightarrow D^{\vee}$ induit des isomorphismes composante par composante, \textit{i.e}. $\lambda_0(D_{\mu})= D_{\mu}\subset D^{\vee}$ et en particulier $\lambda_0(D_{\mu})\subset D_{\mu}^{\perp}$. 

\medskip

Dans cette situation l'espace des filtrations admissibles contient un ouvert de Zariski. Une version à $2$ pentes, qui suffit pour la démonstration du théorème~\ref{theo:EAdm}, est donnée par le lemme suivant.
\begin{lem} \label{lem:zariskcasordin} Si $\mP=\{\mu, 1-\mu\}$ avec $\mu\neq 1/2$ alors
    une filtration $F\subset D_L$ qui vérifie $F\cap D_{\mu}= \{0\}$ et $F\cap D_{1-\mu}= \{0\}$ est admissible. En particulier, l'espace des filtrations admissibles contient un ouvert de Zariski d'intersection non vide avec la grassmannienne lagrangienne.
\end{lem}

\begin{proof} Dans ce cas, tout sous-$\varphi$-module $N$ de $D$ admet une décomposition $N_{\mu}\oplus N_{1-\mu}$ avec $N_{\mu}\subset D_{\mu}$ et $N_{1-\mu}\subset D_{1-\mu}$. La condition d'admissibilité pour une filtration $F\subset D_L$ s'écrit ici 
\[
\dim F\cap N \leq \mu \dim N_{\mu} + (1-\mu) \dim N_{1-\mu}.
\]
    Soit $N\subset D$ un sous-$\varphi$-module.  Alors, l'hypothèse sur $F$ assure que l'on a
    \[
    \dim F\cap N\leq \min (\dim N_{\mu}, \dim N_{1-\mu})
    \]
    et donc l'admissibilité de $F$. 

    \medskip
    
    On utilise alors le lemme~\ref{lem:critlag} pour conclure car $\dim D_{\mu}=\dim D_{1-\mu}=1/2\dim D$. L'ensemble des filtrations considérées est donné par l'intersection de deux ouverts de Zariski non vide de la grassmannienne lagrangienne.
\end{proof}

\begin{rem}
    La conclusion du lemme précédent vaut encore pour un nombre arbitraire de pentes différentes de $1/2$. Pour cela il convient de considérer les parties $I\subset \mP$ de cardinal $1/2\Card \mP$ contenant, pour chaque paire de pentes échangées par la dualité $(\mu, 1-\mu)$, l'une des deux mais pas l'autre et les sous-modules $D_I= \bigoplus_{\mu \in I} D_{\mu}$. Les sous-$\varphi$-modules $D_I$ ainsi obtenus sont lagrangiens pour $\lambda_0$ et un raisonnement analogue à celui présenté dans le lemme conclut. 
\end{rem}

\subsection{Interlude sur les représentations $K$-élémentaires} \label{sub:kelem}

On commence par les définitions essentielles de cette partie et quelques remarques les concernant. On considère dans toute la partie une extension $L/K$ galoisienne de corps de caractéristique $0$ et $G$ un groupe fini.

\begin{defin}
    Soient $V$ et $W$ deux $L$-représentations de $G$. Elles sont dites conjuguées sur $K$ s'il existe un élément $\sigma\in \Gal (L/K)$ et un isomorphisme de $K$-représentations $V\to W$ qui est $\sigma$-semi-linéaire. 
\end{defin}

\begin{rem}  \phantom{v}

\begin{itemize}
    \item[$(i)$] La conjugaison définit une  relation d'équivalence sur les classes d'isomorphie de $L$-représentations de $G$. 
    \item[$(ii)$] Deux $L$-représentations $V$ et $W$ de $G$ sont conjuguées sur $K$ s'il existe un élément $\sigma\in \Gal (L/K)$ et un isomorphisme de $L$-représentations 
    \[
    V\otimes_{L^{\sigma}} L \to W. 
    \]
\end{itemize}
\end{rem}

\begin{defin}
    Une $L$-représentation de $G$ est dite $K$-élémentaire si, pour tout couple de composantes isotypiques $(W,W')$ de $V$, $W'$ est conjuguée sur $K$ à $W$ ou $W^{\vee}$. 
\end{defin}

\begin{rem}  \phantom{v}

\begin{itemize}
    \item[$(i)$] Dans une représentation $K$-élémentaire toutes les composantes isotypiques ont la même dimension.
    \item[$(ii)$] Il suit de la définition que la propriété est encore vraie pour tout couple de sous-représentation irréductibles $(W,W')$ d'une représentation $K$-élémentaire $W$, c'est-à-dire que pour un tel couple $(W,W')$, $W'$ est conjuguée à $W$ ou $W^{\vee}$. 
\end{itemize}
\end{rem}

On commence par quelques lemmes sur les propriétés basiques des représentations $K$-élémentaires. Rappelons que, pour $V$ un $K$-espace vectoriel, le groupe de Galois $\Gal (L/K)$ agit naturellement sur $V\otimes_K L$ à droite par $K$-automorphismes. Les conjugués d'un sous-$L$-espace vectoriel $W\subset V\otimes_K L$ sont les $\sigma(W)$ pour $\sigma\in \Gal(L/K)$.

\begin{lem} \label{lem:repkelemirred}
    Soient $V$ une $K$-représentation irréductible (respectivement isotypique) de $G$. Alors $V\otimes_K L$ est une $L$-représentation qui est $K$-élémentaire et s'écrit comme la somme directe de conjugués de toute sous-représentation irréductible (respectivement isotypique maximale) $W\subset V\otimes_K L$. 
\end{lem}

\begin{proof}
    Considérons d'abord le cas $V$ irréductible. Soit $W\subset V\otimes_K L$ une sous-représentation irréductible. Alors, la somme 
    \[
    \sum\limits_{\sigma\in \Gal (L/K)} \sigma(W)
    \]
    est une sous-représentation de $V\otimes_K L$ qui provient de $K$, c'est-à-dire qu'il existe $W'\subset V$ une $K$-représentation de $G$ telle que 
    \[
    W'\otimes_K L =\sum\limits_{\sigma\in \Gal (L/K)} \sigma(W).
    \]
    Il suit que $W'=V$ par irréductibilité et le résultat voulu en extrayant une somme directe de la somme donnée pour $W'\otimes_K L$. Le cas où $V$ est isotypique suit.
\end{proof}

\begin{lem} \label{lem:changebasekelem}
    Soit $V$ une $L$-représentation de $G$ qui est $K$-élémentaire. Alors la $\overline{K}$-représentation $V\otimes_L \overline{K}$ de $G$ est $K$-élémentaire. 
\end{lem}

\begin{proof}
    On considère la décomposition en composantes isotypiques
    \[
    V\otimes_L \overline{K} = \bigoplus\limits_{i=1}^r W_i.
    \]
    On peut supposer que $r\geq 2$, le cas $r=1$ étant trivial. On considère la décomposition en composantes isotypiques
    \[
    V=\bigoplus_{s=1}^{t} V_s.
    \]
    Il suit du lemme~\ref{lem:repkelemirred} que, pour tout $s\in \{1,\dots, t\}$, $V_s\otimes_L \overline{K}=\bigoplus_{k=1}^{r_s} W_{a_k}$ pour une famille d'indice $\{a_1,\dots, a_{r_s}\}\in \{1,\dots, r_s\}$. De plus, $W_{a_k}$ est conjugué à $W_{a_{k'}}$ sur $L$ pour $k,k'\in \{1,\dots, r_s\}$. En effet, pour une sous-représentation irréductible de $V\otimes_L \overline{K}$, la somme de ses conjugués galoisiens déterminent la sous-représentation de la forme $V_s\otimes_L \overline{K}$ à laquelle elle appartient.  

    \medskip

    On montre maintenant que $V\otimes_L \overline{K}$ est $K$-élémentaire. Soient donc $i,j\in \{1,\dots, r\}$ avec $i\neq j$ et $(W_i,W_j)$ le couple correspondant de composantes isotypiques de $V\otimes_L \overline{K}$. Tout d'abord, s'il existe $s\in \{1,\dots ,t\}$ tel que $W_i$ et $W_j$ sont des sous-représentations de $V_s\otimes_L \overline{K}$ alors il existe un isomorphisme de $L$-représentations $\sigma$-semi-linéaire $W_i\to W_j$ pour un $\sigma\in \Gal(\overline{K}/L)$. C'est, en particulier, un isomorphisme de $K$-représentations et cela conclut dans ce cas.  

    \smallskip
    
    Sinon, il existe $s,s'\in \{1,\dots, t\}$ tels que $W_i\subset V_s\otimes_L \overline{K}$ et $W_j\subset V_{s'}\otimes_L \overline{K}$. Comme $V$ est $K$-élémentaire, il existe un $K$-isomorphisme $\sigma'$-semi-linéaire $V_s\to V_{s'}$ (ou $V_s\to V_{s'}^{\vee}$) avec $\sigma' \in \Gal(L/K)$. Celui-ci se relève en $K$-isomorphisme $\sigma$-semi-linéaire $f\colon V_s\otimes_L \overline{K} \to V_{s'}\otimes_L \overline{K}$ (ou $ V_s\otimes_L \overline{K} \to V_{s'}^ {\vee}\otimes_L \overline{K}=(V_{s'}\otimes_L \overline{K})^{\vee}$) pour un relevé $\sigma \in \Gal (\overline{K}/K)$ de $\sigma'$. L'image par $f$ de $W_i$ est l'un des $W_{s'_k}$ (ou $W_{s'_k}^{\vee}$) et on conclut en composant par l'isomorphisme de $K$-représentation déduit du cas précédent en remplaçant $W_i$ par $f(W_i)$. 
    
\end{proof}

\begin{lem} \label{lem:factorepkelem}
    Soient $V$ une $L$-représentation $K$-élémentaire de $G$ et $W$ une sous-représentation irréductible. Alors le morphisme $G\to GL(V)$ se factorise à travers le morphisme $G\to GL(W)$. 
\end{lem}
\begin{proof}
    L'écriture de $V$ en somme directe de représentations irréductibles $V=\bigoplus_{i=1}^r W_i$ donne une factorisation de $G\to \GL(V)$ par $G\to \prod_{i=1}^r \GL(W_i)$. Il suffit donc de montrer que, pour $i\in \{1,\dots, r\}$ fixé, $G\to \GL(W_i)$ se factorise par $G\to GL(W)$. Si $W$ et $W_i$ sont conjuguées par $T_i\colon W\to W_i$ sur $K$, on considère l'isomorphisme $\GL(W)\to \GL(W_i)$ donné par $f\mapsto T_i \circ f \circ T_i^{-1}$. Si $W$ est conjuguée au dual de $W$ alors on compose de plus avec l'isomorphisme
    \[
    \GL(W)\to \GL(W^{\vee}),~f\mapsto (\prescript{t}{}{f})^{-1}.
    \]
\end{proof}

On va maintenant s'intéresser au lien entre les représentations $K$-élémentaires de $G$ et une classe particulière d'endomorphismes d'espaces vectoriels.

\begin{defin}
    Un endomorphisme $f$ d'un $K$-espace vectoriel de dimension finie $V$ est dit perturbateur pour $V$ si tous les sous-espaces propres de l'endomorphisme $f_{\overline{K}}$ de $V\otimes_K \overline{K}$ sont de dimension au plus $(\dim V)/2$.  
\end{defin}

Lorsque le contexte est clair, pour une $K$-représentation $\rho\colon G\to \GL(V)$ de dimension finie on dira simplement que $h\in G$ est perturbateur pour $V$ si $\rho(h)$ l'est. On a directement le résultat suivant.

\begin{lem}\label{lem:decompkelem}
    Soient $V$ et $W$ deux $L$-représentations de $G$ et $h\in G$ un élément perturbateur pour $V$ et $W$. Alors $h$ est également perturbateur pour les représentations $V\oplus W$, $V^{\vee}$ et toutes les représentations conjuguées de $V$ sur $K$. 
\end{lem}

Le lemme qui suit permet de vérifier l'existence d'éléments perturbateurs pour une représentation de $G$. 

\begin{lem}\label{lem:rep1/2}
    Si $K$ est algébriquement clos et $V$ est une $K$-représentation irréductible de $G$ de dimension $d\geq 2$, alors il existe un élément de $G$ perturbateur pour $V$. 
\end{lem}
\begin{proof}
    On peut supposer $K=\C$. Soit $\chi$ le caractère de $\rho$. Par irréductibilité on a 
    \[
    \sum\limits_{g\in G} |\chi(g)|^2= \Card G.
    \]
    Comme $\chi(e)=\dim V\geq 2$, il existe $g\in G$ avec $|\chi(g)|<1$. On note $\lambda_1,\dots,\lambda_n$ les valeurs propres de $\rho(g)$ et $E_{\lambda_1},\dots E_{\lambda_n}$ les sous-espaces propres associés. Du fait que $\rho(g)$ est diagonalisable on a 
    \[
    \chi(g)=\sum\limits_{i=1}^n \dim E_{\lambda_i} \lambda_i.
    \]
    Comme les $(\lambda_i)_{i\in \{1,\dots n\}}$ sont des racines de l'unité, par inégalité triangulaire on a
    \[
    \dim E_{\lambda_i} \leq |\chi(g)|+ \sum\limits_{j\neq i} \dim E_{\lambda_j}.
    \]
    On en déduit l'inégalité voulue du fait que $\sum_{j\neq i} \dim E_{\lambda_j}= \dim V -\dim E_{\lambda_i}$. 
\end{proof}

La proposition principale sur les représentations $K$-élémentaires pour notre propos est la suivante.

\begin{prop}\label{prop:vpgalois}
    Soit $V$ une $L$-représentation $K$-élémentaire de $G$. Alors, ou bien il existe un élément de $G$ perturbateur pour $V$ ou bien $G$ agit par homothéties.
\end{prop}
\begin{proof}
    D'après la définition d'élément perturbateur et le lemme~\ref{lem:changebasekelem} on peut supposer $L$ algébriquement clos. Supposons d'abord que la dimension commune des sous-représentations irréductibles de $V$ est au moins $2$. Dans ce cas, les lemmes~\ref{lem:decompkelem} et~\ref{lem:rep1/2} montrent qu'il existe un élément de $G$ perturbateur pour $V$. Si, au contraire, cette dimension commune vaut $1$ alors, d'après le lemme~\ref{lem:factorepkelem} le morphisme $G\to \GL(V)$ se factorise par un morphisme $G\to L^{\times}$ et donc a son image cyclique dans $\GL(V)$. Soit $h\in G$ qui engendre cette image. Alors, les sous-espaces propres pour $h$ sont exactement les composantes isotypiques de $V$ et ont donc même dimension par $K$-élémentarité. Il suit que, ou bien $h$ est perturbateur pour $V$ ou bien $h$ n'a qu'un seul sous-espace propre et donc $G$ agit par homothétie, son image dans $\GL(V)$ étant engendrée par $h$. 
\end{proof}

\subsubsection{} Pour finir, on montre le critère annoncé en partie~\ref{sub:supersing} pour vérifier l'hypothèse du théorème~\ref{theo:EAdmss}.

\begin{lem} \label{lem:existlagvpcond}
    Soient $(V,\lambda)$ un espace symplectique sur $K$ de dimension $2g$ et $h\in \Sp(\lambda)$ un automorphisme symplectique tel que $h_{\overline{K}}$ est diagonalisable et perturbateur pour $V$. Alors il existe un lagrangien $F$ tel que $\dim F\cap h(F) \leq 1$. 
\end{lem}
\begin{proof}
    Tout d'abord, du fait que la grassmannienne lagrangienne est rationnelle, il suffit de démontrer l'existence d'un tel $F$ dans $V_{\overline{K}}$.

    \medskip

    Il existe une base symplectique $(x_i,y_i)_{1\leq i \leq g}$ de $V_{\overline{K}}$ formées de vecteurs propres pour $h$ et qui donne une décomposition de l'espace en somme orthogonale stable par $h$
    \[
    V_{\overline{K}}= \bigoplus\limits_{i=1}^g V(x_i,y_i)
    \]
    où $V(x_i,y_i)$ est un espace symplectique de dimension $2$ standard, c'est-à-dire qui vérifie $\lambda(x_i,x_i)=0$, $\lambda(x_i,y_i)=1$ et $\lambda(y_i,y_i)=0$. On va construire $F$ comme somme directe orthogonale de lagrangiens de sous-espaces symplectiques stables par $h$ de la forme $V(x_i,y_i)$ ou $V(x_i,y_i)\oplus V(x_j,y_j)$ avec $i\neq j$. 

    \medskip

    Soient $\mu_1,\dots, \mu_n$ les valeurs propres de $h$ classées par multiplicité décroissante. On montre d'abord que, si $g\geq 4$ alors il existe un sous-espace symplectique $E\subset V_{\overline{K}}$ de la forme annoncée ainsi qu'un lagrangien $F\subset E$ qui vérifie $h(F)\cap F=\{0\}$ de manière à ce que $E^{\perp}$ et $h_{|E^{\perp}}$ vérifient encore les hypothèses de l'énoncé.

    \medskip

    Tout d'abord, comme les $x_i$ et $y_i$ sont des vecteurs propres de $h$ on peut noter $V(x_i,y_i)$ par $V(\mu_k,\mu_l)$ pour deux valeurs propres (non nécessairement distinctes) $\mu_k$ et $\mu_l$ de $h$ pour signifier que $h(x_i)=\mu_k x_i$ et $h(y_i)=\mu_l y_i$. On va maintenant montrer que l'on peut toujours choisir le couple $(E,F)$ de l'une des deux formes suivantes 
    \begin{itemize}
        \item $E= V(\mu_1,\mu_l)$ avec $\mu_1\neq \mu_l$ et $F$ est donné par un élément de l'intersection des ouverts donnés par le lemme~\ref{lem:critlag} avec $M$ chacune des droites propres pour $h$ de $E$;
        \item $E=V(\mu_1, \mu_1)\oplus V(\mu_l,\mu_m)=V(x_i,y_i)\oplus V(x_j,y_j)$ avec $\mu_1\neq \mu_l$ et $\mu_1\neq \mu_m$ et $F$ est le lagrangien engendré par $x_i+x_j$ et $y_i-y_j$. Dans ce cas, il est clair que $F$ est lagrangien et il faut vérifier que $h(F)\cap F=\{0\}$. Un élément de cette intersection vérifie l'égalité
        \[
        ax_i+ax_j+by_i-by_j=c\mu_1 x_i+c\mu_l x_j+d\mu_1 y_i-d\mu_m y_j
        \]
        pour $a,b,c,d\in \overline{K}$. On en déduit directement $a=b=c=d=0$ du fait que $\mu_1\neq \mu_l$ et $\mu_1\neq \mu_m$. 
    \end{itemize}
    On procède par cas, suivant la multiplicité des valeurs propres de $h$. Si $h$ n'admet que $2$ valeurs propres (toutes deux de multiplicité $g$) ou si toutes les multiplicités des valeurs propres de $h$, sauf au plus une, sont plus petites que $g-2$ alors tout choix de $E$ d'une des deux formes convient et nécessairement au moins l'une des deux existe. Il faut donc traiter le cas où $h$ a une unique valeur propre de multiplicité $g$, une de multiplicité $g-1$ et une de multiplicité $1$ et celui où $h$ admet au moins deux valeurs propres de multiplicité $g-1$. 

    \medskip 

    Dans le premier cas, la valeur propre de multiplicité $1$, $\mu_3$ apparaît dans un $V(x_i,y_i)$ de la forme $V(\mu_k,\mu_3)$. Si $k=1$ alors ce choix pour $E$ convient. Si $k=2$ alors il existe nécessairement un $V(x_i,y_i)$ de la forme $V(\mu_1,\mu_1)$ et $E=V(\mu_1,\mu_1)\oplus V(\mu_2,\mu_3)$ convient. Dans le deuxième cas, comme il n'y a pas de valeurs propres de multiplicité $g$, s'il existe sous-espace de la première forme il convient. S'il n'en existe pas, on peut considérer $E$ de la forme $V(\mu_1,\mu_1)\oplus V(\mu_2,\mu_2)$. Si celui-ci ne convient pas, c'est que la valeur propre $\mu_3$ est de multiplicité $g-1$, donc que $g$ vérifie l'inégalité $2g\geq 3(g-1)$ et par suite $g = 3$.

    \medskip

    Par une récurrence immédiate on est donc ramené à traiter les cas $g\leq 3$. Si $g=1$ ou $g=2$ en utilisant à nouveau la construction précédente on voit qu'il existe $F$ lagrangien qui vérifie $F\cap h(F)=\{0\}$. Si $g=3$, le seul cas restant est celui où $h$ admet trois valeurs propres de multiplicité $2$. Dans ce cas, on obtient $F$ comme somme orthogonale de $F_1$ sur $V(\mu_1,\mu_1)\oplus V(\mu_2,\mu_2)$ construit comme précédemment et de n'importe quel lagrangien $F_2$ de $V(\mu_3,\mu_3)$. La somme orthogonale $F_1\oplus F_2$ est lagrangienne et vérifie $\dim F\cap h(F)\leq 1$ par construction. 

    \medskip

    Dans tous les cas, on a construit $F$ qui vérifie les hypothèses de l'énoncé. 
    
\end{proof}

\subsection{Conclusion}

\subsubsection{} Le but de cette partie est de montrer que pour une extension galoisienne $L/K$ fixée comme au paragraphe~\ref{subsub:galoisextexist}, il existe toujours des $L$-filtrations vérifiant les propriétés (\ref{enum:lag}), (\ref{enum:adm}) et (\ref{enum:stabgal}), c'est-à-dire de démontrer le théorème~\ref{theo:EAdm}. D'après la proposition~\ref{prop:reductab} il suffit de traiter le cas où $D$ est un $\varphi$-module abélien. On va montrer que le problème se réduit encore suivant une décomposition orthogonale pour $\lambda_0$, $G$-stable et $\varphi$-stable. On commence par établir le lien structurel entre $\varphi$-modules et représentations $K$-élémentaires. 

\begin{lem}\label{lem:structphiG}
    Soient $D$ un $\varphi$-module polarisé muni d'une action d'un groupe fini $G$ et $V\subset D$ un sous-espace vectoriel. Alors, les assertions suivantes sont vérifiées.
    \begin{enumerate}
        \item Si $V$ est $\varphi$-stable (respectivement $G$-stable) alors $V^{\perp}$ est $\varphi$-stable (respectivement $G$-stable).
        \item \label{enum:phietG} Si $V$ est $G$-stable alors $\varphi(V)$ est $G$-stable et est conjugué à $V$ par $\sigma$. 
        \item Si $V$ est une composante isotypique de $D$ pour l'action de $G$ alors le plus petit sous-$\varphi$-module de $D$ contenant $V$ est une représentation $K$-élémentaire de $G$.
    \end{enumerate}
\end{lem}

\begin{proof}
    Les assertions sur la stabilité s'obtiennent par un calcul direct. Pour $(\ref{enum:phietG})$, il suffit de remarquer que $\varphi$ est un isomorphisme de $K$-représentation $\sigma$-semi-linéaire entre $V$ et $\varphi(V)$. Le dernier point se déduit directement du deuxième en remarquant qu'il existe $n\in \N\setminus \{0\}$ tel que $\varphi^{n}(V)=V$. 
\end{proof}

On continue avec l'existence d'une décomposition qui permettra d'appliquer les lemmes précédents et les résultats sur les représentations $K$-élémentaires pour démontrer le théorème~\ref{theo:EAdm}.

\begin{prop}\label{prop:decomp}
    Soient $D$ un $\varphi$-module abélien polarisé et $G\subset \Aut(D,\lambda_0)$ un groupe fini.

    Alors il existe une décomposition 
    \[
    D= \bigoplus\limits_{i=1}^r D_i
    \]
    telle que la somme est orthogonale pour $\lambda_0$, les $D_i$ sont des sous-$\varphi$-modules abéliens de $D$ dont l'ensemble des pentes est de cardinal au plus $2$ et sont des sous-représentations $\Q_p$-élémentaires de $G$ pour $1\leq i\leq r$.  
\end{prop}
\begin{proof}
    On note $\mP=\{\mu_1,\dots, \mu_a\}$ l'ensemble des pentes de $D$ rangées par ordre croissant. On a une décomposition en composantes isoclines
    \[
    D=\bigoplus\limits_{\mu \in \mP} D_{\mu}.
    \]
    On pose $D'_i=D_{\mu_i}+D_{1-\mu_i}$ pour $1\leq i\leq \big\lceil \frac{a}{2} \big\rceil$. Alors l'écriture $D=\bigoplus_{i=1}^{\big\lceil \frac{a}{2} \big\rceil} D'_i$ est une somme orthogonale de $\varphi$-modules $G$-stables qui vérifient la condition que les pentes de chaque composante sont données par un ensemble de cardinal au plus $2$. On peut donc supposer dans la suite que $\mP=\{\mu,1-\mu\}$ pour un $\mu\in \Q\cap [0,1]$. 

    \medskip

    On procède alors par récurrence suivant la dimension de $D$. Si $D$ est nul, il n'y a rien à faire. Sinon, on considère la décomposition 
    \[
    D=\bigoplus\limits_{i=1}^r V_i
    \]
    de $D$ en composantes isotypiques pour l'action de $G$. On note $A$ la somme des conjuguées de $V_1$ ainsi que leurs duales et on note $B$ le supplémentaire $G$-stable de $A$, qui est unique par construction. Les deux sous-espaces $A$ et $B$ de $D$ sont en particulier $\varphi$-stables d'après le lemme~\ref{lem:structphiG}. Par construction $A$ est $\Q_p$-élémentaire en tant que représentation de $G$ et la restriction de $\lambda_0$ à $A$ induit une forme bilinéaire non dégénérée. Il suit que $A^{\perp}\cap A =\{0\}$. L'orthogonal $A^{\perp}$ de $A$ est un supplémentaire $G$-stable de $A$ et alors $A^{\perp}=B$. L'hypothèse de récurrence s'applique pour $B$ ce qui conclut.
 \end{proof}

\subsubsection{} Pour les $\varphi$-modules abéliens isoclines on peut désormais conclure.

\begin{lem}\label{lem:existQpelem}
    Soit $D$ un $\varphi$-module abélien polarisé isocline, nécessairement de pente $1/2$, et $G\subset \Aut (D,\lambda_0)$ un groupe fini qui fait de $D$ une $G$-représentation $\Q_p$-élémentaire. Alors il existe une filtration admissible sur $D_L$ qui vérifie les propriétés~(\ref{enum:lag}), (\ref{enum:adm}) et (\ref{enum:stabgal}) du théorème~\ref{theo:EAdm}.
\end{lem}

\begin{proof}
    D'après la proposition~\ref{prop:vpgalois} ou bien $G$ agit par homothétie ou bien il existe un élément de $G$ perturbateur pour $D$. Dans le cas où $G$ n'agit pas sur $D$ par homothéties on obtient la filtration grâce au théorème~\ref{theo:EAdmss} qui s'applique d'après le lemme~\ref{lem:existlagvpcond}. Lorsque $G$ agit par homothéties, donc trivialement sur les sous-espaces vectoriels, l'existence d'une filtration admissible définie sur $K_0$ provient du corollaire~\ref{cor:supersingdense}.
\end{proof}

Finalement, on montre la stabilité par décomposition du problème. On précise un peu en quoi cela consiste pour la propriété~(\ref{enum:stabgal}) de stabilité galoisienne. Soient, comme précédemment, un $\varphi$-module abélien polarisé $D$ muni d'une action d'un groupe de ramification $G$ et $L/K$ une extension galoisienne de groupe $G$. Supposons de plus que la décomposition fournie par la proposition~\ref{prop:decomp} s'écrive $D=D_1\oplus D_2$. Pour obtenir une filtration sur $D_L$ avec la propriété~(\ref{enum:stabgal}) on voudrait s'intéresser aux filtrations sur $D_1$ ayant cette même propriété pour l'action induite. Il faut alors considérer le sous-corps $L_1$ de $L$ stabilisé par le noyau de l'action induite sur $D_1$ et considérer la propriété~(\ref{enum:stabgal}) pour le groupe quotient, groupe de Galois de $L_1$ sur $K$. Pour remonter, il suffit alors de constater que si $\sigma\in \Gal(L/K)$ est un relevé de $\tau\in \Gal (L_1/K)$ alors
\[
\sigma\cdot_{lin} F= \tau\cdot_{lin} F \text{ et } \sigma\cdot_{gal} F =\tau\cdot_{gal} F
\]
pour tout sous-espace vectoriel $F\subset (D_1)_{L_1}$ vu comme sous-espace de $D_L$. 

\begin{lem}\label{lem:decomp}
Soient $D$ un $\varphi$-module abélien polarisé et $G\subset \Aut (D, \lambda_0)$ un groupe fini de ramification et une extension $L/K$ galoisienne de groupe $G$. On suppose que $D$ admet une décomposition $D=D_1\oplus D_2$ qui soit orthogonale pour $\lambda$, stable par $G$ et par $\varphi$ et qu'il existe des filtrations $F_1$ et $F_2$ sur $D_1$ et $D_2$ qui vérifient les propriétés (\ref{enum:lag}), (\ref{enum:adm}) et (\ref{enum:stabgal}), la dernière propriété étant pour l'action induite de $G$ sur $D_1$ ou $D_2$ et le sous-corps de $L$ correspondant. Alors la filtration $F=(F_1)_L\oplus (F_2)_L$ sur $D_L$ vérifie les propriétés (\ref{enum:lag}), (\ref{enum:adm}) et (\ref{enum:stabgal}).
\end{lem}
\begin{proof}
    La stabilité de la propriété~(\ref{enum:lag}) par somme orthogonale est triviale. On vérifie les deux autres. Pour la propriété~(\ref{enum:stabgal}), il faut montrer que pour tout $h\in G$ on a l'égalité
    \[
    h\cdot_{gal}(F)= h\cdot_{lin} (F).
    \]
    Soit $h\in G$. Comme les sous-$\varphi$-modules de $D$ sont définis sur $K_0$, l'action galoisienne de $h$ respecte la décomposition $D_1\oplus D_2$ et on a
    \begin{align*}
        h\cdot_{gal}(F)&= h\cdot_{gal}(F_1) \oplus h\cdot_{gal} (F_2) \\
        &= h\cdot_{lin}(F_1)\oplus h\cdot_{lin}(F_2).
    \end{align*}
    La deuxième égalité est donnée par l'hypothèse sur $F_1$ et $F_2$. Finalement, par linéarité on a $h\cdot_{lin}(F_1)\oplus h\cdot_{lin}(F_2)= h\cdot_{lin}(F)$.

    \medskip

    Il reste à vérifier la propriété d'admissibilité. Soit $N\subset D$ un sous-$\varphi$-module. La seconde projection $p_2\colon D\to D_2$ donne une suite exacte de $\varphi$-modules 
    \[
    \begin{tikzcd}
        0\arrow[r] & N\cap D_1 \arrow[r] & N \arrow[r] & p_2(N) \arrow[r] & 0
    \end{tikzcd}
    \]
    qui se casse suivant les composantes isoclines de $D$: pour tout $\mu\in [0,1]\cap \Q$, on a
    \[
    \begin{tikzcd}
        0\arrow[r] & (N\cap D_1)_{\mu} \arrow[r] & N_{\mu} \arrow[r] & p_2(N)_{\mu} \arrow[r] & 0.
    \end{tikzcd}
    \]

D'un autre côté, on dispose d'une suite exacte de $L$-espaces vectoriels
    \[
    \begin{tikzcd}
        0\arrow[r] & N_L\cap F_1 \arrow[r] & N_L\cap F \arrow[r] & p_2(N)_L\cap F_2    \end{tikzcd}
    \]
d'où une inégalité
\[
\dim F\cap N_L \leq \dim N_L\cap F_1+ \dim p_2(N)_L \cap F_2.
\]
L'admissibilité de $F_1$ et $F_2$ assure alors que
\[
\dim F\cap N_L \leq \sum\limits_{\mu} \mu\dim (N\cap D_1)_{\mu}+ \sum\limits_{\mu} \mu\dim p_2(N)_{\mu}= \sum\limits_{\mu} \mu \dim N_{\mu}
\]
ce qui est la condition d'admissibilité de $F$ pour $N$. 
\end{proof}

On peut désormais démontrer le théorème~\ref{theo:EAdm}, démonstration qui consiste à assembler les résultats de la partie~\ref{sec:phimod}.

\begin{proof}
    On considère la décomposition
    \[
    D=\bigoplus\limits_{i=1}^n D_i
    \]
    donnée par la proposition~\ref{prop:decomp}. D'après le lemme~\ref{lem:decomp}, il suffit de démontrer l'existence de filtrations $F_1,\dots, F_n$ sur chacun des $\varphi$-modules $D_1,\dots, D_n$ qui vérifient les propriétés (\ref{enum:lag}), (\ref{enum:adm}) et (\ref{enum:stabgal}) pour les actions induites de $G$. Pour $i\in \{1,\dots, n\}$ tel que l'ensemble des pentes de $D_i$ est de cardinal $2$, l'existence d'une telle filtration $F_i$ s'obtient grâce au lemme~\ref{lem:zariskcasordin}. Pour les indices restants, la filtration $F_i$ s'obtient par le lemme~\ref{lem:existQpelem}.
\end{proof}

\section{Réalisation de groupes finis comme groupes de monodromie finie par déformation} \label{sec:constru}

Le but de cette partie est de donner une méthode de construction de variétés abéliennes avec monodromie finie prescrite, dont le résultat formel est l'énoncé \ref{theo:resumconstruction}. Plus précisément, on montre que pour une donnée $(A_0,\lambda_0, G)$ où $A_0$ est une variété semi-abélienne sur un corps fini $k$ de caractéristique $p$, $\lambda_0$ une polarisation de $A_0$ et $G\subset \Aut (A_0, \lambda_0)$ un groupe fini de ramification, il existe une variété abélienne $B$, relèvement à isogénie géométrique près de $A_0$, sur un corps $p$-adique $K$ telle que le groupe de monodromie finie de $B$ soit $G$. 

\subsection{Quelques rappels de théorie de Hodge $p$-adique entière} \label{sub:rappelsHodge}

Le but de cette partie est de décrire les différentes catégories semi-linéaires et foncteurs qui interviennent dans la construction ainsi que d'établir le lien avec la partie précédente.

\subsubsection{} Pour un corps fini $k$ et $K_0$ l'extension maximale non ramifiée de  $\operatorname{Frac}(W(k))$ on a introduit en partie \ref{sec:phimod} les catégories de Fontaine des $\varphi$-modules $\operatorname{MF}^{\varphi}$ et celle des $\varphi$-modules filtrés sur $K$ et admissibles $\operatorname{MF}^{\varphi}_K$ pour une extension finie $K/K_0$. 

\medskip

D'après le théorème~A p. 140 de \cite{FO22} on dispose d'une équivalence de catégorie entre $\operatorname{MF}^{\varphi,N}_K$ et la catégorie $\operatorname{Rep}^{st}_{\Q_p}(G_K)$ des représentations $p$-adiques de $G_K$ qui sont semi-stables. Cette équivalence est donnée par deux foncteurs, définis sur les objets comme suit
\[
\begin{array}{ccccccc}
  D_{st}\colon  \operatorname{Rep}^{st}_{\Q_p}(G_K) & \longrightarrow & \operatorname{MF}^{\varphi,N}_K & \text{et} & V_{st}\colon \operatorname{MF}^{\varphi,N}_K  & \longrightarrow & \operatorname{Rep}^{st}_{\Q_p}(G_K)  \\
      V & \longmapsto & (V\otimes_{\Q_p} B_{st})^{G_K} & & D & \longmapsto & D\otimes_{K_0} B_{st} \\  
\end{array}
\]
où $B_{st}$ est un anneau de périodes défini en partie~6.1 de \cite{FO22}.
Par cette équivalence les représentation cristallines, dont on note la catégorie $\operatorname{Rep}^{\mathrm{cris}}_{\Q_p}(G_K)$, sont équivalentes à la sous-catégorie pleine des $\varphi$-modules, c'est-à-dire les objets qui vérifient $N=0$. De plus, les poids de Hodge-Tate d'une telle représentation $V_{st}(D)$ pour un $\varphi,N$-module $D$ sont donnés par les entiers $i\in \Z$ tels que $F_{i+1}\neq F_i$ au niveau de la filtration $F=\{F_i\}_{i\in \Z}$ de $D_K$. Les $\varphi$-modules étudiés en partie~\ref{sec:phimod} correspondent à des représentations cristallines à poids de Hodge-Tate $0,1$ par construction.

\subsubsection{} On note $\mfS$ l'anneau des séries formelles $W(k)[[u]]$ et $E(u)$ un polynôme d'Eisenstein pour l'extension $K/K_0$. Il est naturellement équipé d'un morphisme de Frobenius $\varphi\colon \sum_{n\geq 0} a_n u^n \mapsto \sum_{n\geq 0} \sigma(a_n) u^{pn}$ où $\sigma$ est le morphisme de Frobenius canonique sur $W(k)$ donné par l'élévation à la puissance $p$ sur les coordonnées. La catégorie des modules de Breuil-Kisin, introduite dans \cite{Kis06}, est une sous-catégorie notée $\operatorname{BT}_{\mfS}^{\varphi}$ de la catégorie $\operatorname{Mod}_{\mfS}^{\varphi}$ des $\mfS$-modules libres de rang fini équipés d'un morphisme injectif $f$ semi-linéaire par rapport à $\varphi$ dont le conoyau est    tué par une puissance de $E(u)$. La sous-catégorie $\operatorname{BT}_{\mfS}^{\varphi}$ des modules de Breuil-Kisin est celle dont les objets ont leur Frobenius $f$ de conoyau tué par $E(u)$.  

\medskip

Kisin montre que la catégorie d'isogénie $\operatorname{Mod}_{\mfS}^{\varphi}\otimes \Q_p$ des modules de Kisin est équivalente à la catégorie $\operatorname{MF}^{\varphi,0,1}_K$ des $\varphi$-modules filtrés sur $K$ de Fontaine dont la filtration à un cran introduite en partie~\ref{sec:phimod} -- voir \cite{Kis06} proposition~2.2.2 -- par la donnée d'un foncteur
\[
\oD \colon \operatorname{BT}_{\mathfrak{S}}^{\varphi} \longrightarrow \mathrm{MF}_{K}^{\varphi,0,1}
\]
Cette équivalence passe par l'extension des scalaires des modules de Kisin à l'anneau des fonctions holomorphes $\OO$ sur le disque unité dans la théorie analytique rigide de Tate. L'anneau $\OO$ se décrit comme le sous-anneau des séries formelles à coefficients dans $K_0$ qui convergent en les $\overline{K}$-points du disque unité. Il construit ensuite deux foncteurs, comparables aux foncteurs $D_{st}$ et $V_{st}$ de Fontaine, 
\[
D_{\mfS}\colon \operatorname{Rep}^{\mathrm{cris~0,1}}_{\Z_p}(G_K) \to \operatorname{BT}_{\mfS}^{\varphi}, ~ V_{\mfS} \colon \operatorname{BT}_{\mfS}^{\varphi} \to \operatorname{Rep}^{\mathrm{cris~0,1}}_{\Z_p}(G_K)
\]
où la catégorie $\operatorname{Rep}^{\mathrm{cris~0,1}}_{\Z_p}(G_K)$ est celle des réseaux dans les représentations cristallines à poids de Hodge-Tate $0,1$.
Ceux-ci sont compatibles aux foncteurs de Fontaine au sens où ils referment et font commuter le diagramme suivant par la ligne en pointillés 
\[
\begin{tikzcd}
    \operatorname{Rep}^{\mathrm{cris~0,1}}_{\Q_p}(G_K) \arrow[r] & \arrow[l]  \operatorname{MF}^{\varphi,0,1}_K \\
    \operatorname{Rep}^{\mathrm{cris~0,1}}_{\Z_p}(G_K) \arrow[u] \arrow[r, dashed] & \arrow[l, dashed] \operatorname{BT}_{\mfS}^{\varphi} \arrow[u]
\end{tikzcd}
\]

\subsubsection{} On rappelle la construction d'un foncteur fidèle $\mD$ de la catégorie $\operatorname{BT}_{\mathfrak{S}}^{\varphi}$ dans celle des réseaux dans les $\varphi$-modules filtrés de Fontaine notée $\mathrm{MF}_{K,W}^{\varphi}$. Une telle construction est par exemple donnée dans la preuve du théorème 2.3 de \cite{Liu12}. Précisément, la catégorie $\mathrm{MF}_{K,W}^{\varphi}$ est celle des $\varphi$-modules $D$ filtrés sur $K$ munis d'un $W(k)$-réseau $W$ tel que $W\otimes_{W(k)} K_0=D$. On note les objets de cette catégorie par un couple $(D,W)$. On définit alors le foncteur $\mD\colon \mathrm{BT}_{\mfS}^{\varphi} \rightarrow \mathrm{MF}_{K,W}^{\varphi}$ de la manière suivante. Sur les objets, à un module de Kisin $\mfM$ on associe le couple $(\oD( \mfM) , W)$ où $W=\mfM/\mfM\cap u\mathrm{M}$ avec $M=\mfM\otimes_{\mfS} \OO$ -- l'anneau $\OO$ étant l'algèbre de Tate introduite précédemment -- et $\oD$ est le foncteur de Kisin rappelé plus haut. Il faut vérifier l'égalité $W\otimes K_0=\oD(\mfM)$ et que ce dernier est stable par $\varphi$. La stabilité de $\mfM$ par $\varphi$ induit celle de $W$. Pour l'autre condition on considère le morphisme surjectif $\psi\colon\mfM/u\mfM\rightarrow W$ et son noyau $\Ker \psi$. Comme $\Ker \psi$ est un sous-$W(k)$-module du module libre de rang fini $\mfM/u\mfM$ il est libre de rang fini et on dispose d'une base adaptée $(\overline{w_1},\dots, \overline{w_n})$. Alors, si $\Ker \psi \neq 0$, la famille $\{p^{n_i} w_i\}_{1\leq i \leq r}$ est une base de $\Ker \psi$ avec $r\geq 1$. Si $(e_1,\dots,e_n)$ une base de $\mfM$ alors les $e_i\mod u$ donnent une base $\{\overline{e_i}\}$ de $\mfM/u\mfM$. Soit une telle base $(e_1,\dots, e_n)$ de $\mfM$. On a une matrice $\overline{P}$ de changement de base de la base $\{\overline{e_i}\}$ à la base $\{\overline{w_i}\}$. Un relevé de $P$ comme matrice à coefficients dans $\mfS$ vérifie $\det P \mod u = \det \overline{P}  \in W(k)^{\times}$ et donc $P$ est inversible. Il suit que $P$ induit un changement de base de la base $\{e_i\}$ à une base $\{w_i\}$ de $\mfM$ avec $w_i\mod u=\overline{w}_i$. Alors $p^{n_1}w_1 \in u\mathrm{M}\cap \mfM$ donc il existe $y\in \mathrm{M}$ tel que $p^{n_1}w_1=uy$. Par ailleurs $y\in \mfM\otimes \OO =\{ \alpha x \mid \alpha \in \OO,~x\in \mfM\}= \mathrm{M}$ ce qui assure que les $\{w_i\}$ forment une base de $\mathrm{M}$. On en déduit $y=\frac{p^{n_1}}{u}w_1$ ce qui est absurde. Le morphisme $\psi$ est donc un isomorphisme et $W$ vérifie $W\otimes_{W(k)} K_0= D$. 

\medskip

Pour les morphismes, à $f\colon \mfM_1\rightarrow \mfM_2$ on associe $\overline{f}=f\mod u$, le morphisme donné par le foncteur $\oD$ de Kisin. On vérifie sans difficulté $\overline{f}(\mfM_1/\mfM_1\cap u\mathrm{M}_1)\subset \mfM_2/\mfM_2\cap u\mathrm{M}_2$. Il suit directement que, comme le foncteur $\oD$ est fidèle, le foncteur $\mD$ l'est aussi.

\subsubsection{} \label{subsub:liensemiabphimod} Soit $A_0$ une variété semi-abélienne sur un corps fini $k$ et polarisée au sens de la définition~\ref{def:vsapol}. Quitte à étendre le corps de base, on peut supposer que les tores en jeu sont déployés. Alors le groupe $p$-divisible $A_0[p^{\infty}]$ fournit un module de Dieudonné $M(A_0)$ sur l'anneau des vecteurs de Witt $W(k)$ -- voir \cite{Fo77} pour cette construction. Par le changement de base à $K_0$ on obtient $D=M(A_0)\otimes_{W(k)} K_0$ un $\varphi$-module issu de $A_0$. Celui-ci admet une décomposition isocline
\[
D= \bigoplus_{\mu\in \Q} D_{\mu}
\]
dont les pentes $\mu\in \Q$ avec $D_{\mu}$ non trivial sont dans $[0,1]$. Par fonctorialité, le $\varphi$-module $D$ s'inscrit dans une suite exacte
\[
0\to D_{T_0} \to D \to D_{B_0} \to 0.
\]
et la polarisation de $A_0$ donne un isomorphisme $\lambda$ de $D$ sur un $\varphi$-module $D^t$ qui induit un isomorphisme de $D_{B_0}$ sur son dual $D_{B_0}^{\vee}$, autrement dit $D_{B_0}$ est un $\phi$-module autodual.

\smallskip

Il suit que $D$ est un $\varphi$-module semi-abélien polarisé auquel le théorème~\ref{theo:EAdm} s'applique. En particulier, pour toute donnée d'un groupe fini $G\subset \Aut (A_0, \lambda_0)$ de ramification on obtient une extension $L/K$ galoisienne de groupe $G$ avec une filtration admissible $F$ sur $D_L$ qui munit $(D,F)$ d'une donnée de descente pour l'extension $L/K$. On clarifie ce dernier point. Si $F\subset D_L$ est une filtration alors le conjugué galoisien du $\varphi$-module filtré $(D,F)$ par un élément $h\in G$ est simplement $(D,h\cdot_{\operatorname{gal}}F)$. Une donnée de descente pour $(D,F)$ correspond alors à une collection d'isomorphismes $(f_h)_{h\in G}$ de $\varphi$-modules filtrés avec $f_h\colon (D,F) \rightarrow (D,h\cdot_{\operatorname{gal}}F)$ vérifiant la condition de cocycle $f_{gh}=g\cdot_{\operatorname{gal}} f_h \circ f_g$ pour $g,h\in G$. 

\smallskip

La condition de cocycle se vérifie alors facilement de la propriété~(\ref{enum:stabgal}) de la filtration $F$ fournie par le théorème~\ref{theo:EAdm} du fait que les morphismes de $\varphi$-modules filtrés commutent à l'action galoisienne. En effet, si $f\colon (D,F)\to (D,H)$ est un tel morphisme alors $h\cdot_{\mathrm{gal}} f$ est simplement $f$ vu comme morphisme entre les $\varphi$-modules filtrés $(D,h\cdot_{\operatorname{gal}} F)$ et $(D, h\cdot_{\operatorname{gal}} H)$, ce qui est bien défini du fait que $f(h\cdot_{\operatorname{gal}} F)= h\cdot_{\operatorname{gal}}(f(F))\subset h\cdot_{\operatorname{gal}} H$ par définition. 

\smallskip

Remarquons finalement que les données de descente qui nous intéressent redonnent l'action du groupe opposé $G^{\mathrm{op}}$ de $G$ sur $D$ par oubli de la filtration, c'est-à-dire que $f_h$ coïncide avec $M(h)$ sur $D$ pour $h\in G$ où $M$ est le foncteur \og{}module de Dieudonné \fg{}.

\medskip

La situation qui nous intéresse en partie~\ref{sub:deformHodge} est décrite par le diagramme commutatif de catégories suivant, où les flèches verticales correspondent au passage \og{} à isogénie près \fg{},

\[
\begin{tikzcd} \label{diag:categories}
    \Rep_{\Q_p}^{\mathrm{cris},0,1}(G_L) \arrow[r] & \operatorname{MF}^{\varphi}_L\simeq\operatorname{BT}_{\mfS}^{\varphi}\otimes \Q &  \\
    \Rep_{\Z_p}^{\mathrm{cris},0,1}(G_L) \arrow[r] \arrow[u] & p-\mathrm{div}/\OO_L\simeq \operatorname{BT}_{\mfS}^{\varphi} \arrow[u] \arrow[r] &  \operatorname{MF}^{\varphi}_{L,W} \arrow[ul]\\
\end{tikzcd}
\]

Tout ceci est de plus compatible avec les passages aux fibres spéciales. En particulier, si $\Gamma$ est un groupe $p$-divisible sur $\OO_L$ alors le $\varphi$-module induit par le module de Dieudonné de sa fibre spéciale est déterminé par la représentation $p$-adique issue de son module de Tate.

\subsection{Déformation de groupes $p$-divisibles issus de variétés semi-abéliennes} \label{sub:deformHodge}

\subsubsection{} Dans cette partie, on démontre un théorème de relèvement en caractéristique $0$ d'une donnée $(A_0, \lambda_0, G)$. En vue d'énoncer le théorème on rappelle quelques faits sur les schémas semi-abéliens sur les anneaux d'entiers de corps $p$-adiques suivant les chapitres 2 et 3 de \cite{FC90} et sur les corps finis ainsi que leurs polarisations. Dans ce texte, pour $L$ un corps $p$-adique, un schéma semi-abélien sur $\OO_L$ est une extension 
\[
\begin{tikzcd}
    0 \arrow[r] & T \arrow[r] & A \arrow[r] & B \arrow[r] & 0
\end{tikzcd}
\]
où $T$ est un tore sur $\OO_L$ et $B$ une variété abélienne sur $\OO_L$. Cette extension correspond à la donnée d'un morphisme $c\colon \underline{X}(T) \to B^{\vee}(\OO_L)$.  

\medskip

\begin{defin}
    Un morphisme sur $\OO_L$ de schémas semi-abéliens $\lambda \colon A\to A^t$ est une polarisation s'il respecte les suites exactes, induit une polarisation $\lambda_B$ entre les quotients abéliens de $A$ et $A^t$, en particulier $B^t=B^{\vee}$ et s'il induit une isogénie $\lambda_T$ sur les parties toriques. Le degré d'une polarisation $\lambda$ d'un schéma semi-abélien est $m^2\cdot \deg \lambda_B$ où $m$ est l'ordre de $\Ker 
    \lambda_T$.  Un schéma semi-abélien sur $\OO_L$ muni d'une polarisation est dit polarisé.
\end{defin}

Un schéma semi-abélien polarisé correspond aux données fournies par les points $(1)$ à $(5)$ de la définition d'une donnée de dégénérescence, objet de la catégorie $\mathrm{DD}_{\mathrm{pol}}$ introduite dans \cite{FC90} p.57--58.

\medskip

Il est clair que la fibre spéciale d'un schéma semi-abélien polarisé sur $\OO_L$ est une variété semi-abélienne polarisée sur $k$. 

\begin{theo} \label{theo:deform}
    Soit $(A_0, \lambda_0)$ une variété semi-abélienne polarisée sur un corps fini $k$. Soit $G\subset \Aut (A_0, \lambda_0)$ un groupe fini de ramification en $p$. Pour toute extension de corps $L/K$ galoisienne de groupe $G$ et de corps résiduel $\overline{k}$ il existe un schéma semi-abélien polarisé $A$ sur $\OO_L$ muni d'une donnée de descente $(f_{\sigma})_{\sigma\in G}$ qui induit une injection $\iota \colon G \hookrightarrow \Aut A_{\overline{k}} \subset \End A_{\overline{k}} \otimes \Q$ par passage à la fibre spéciale où la variété semi-abélienne $A_{\overline{k}}$ est isogène à $(A_0)_{\overline{k}}$ et l'injection $\iota$ correspond à la composée $G\subset \Aut (A_0)_{\overline{k}}\subset \End (A_0)_{\overline{k}}\otimes \Q= \End B_{\overline{k}} \otimes \Q$. De plus, si $\lambda_{B_0}$ est de degré $m$ la polarisation $\lambda$ sur $A$ est telle que la polarisation induite sur le quotient abélien $B$ de $A$ est de degré $mp^n$ pour un $n\geq 0$ et si $\lambda_{T_0}$ est de degré $m'$, l'isogénie induite sur la partie torique $T$ de $A$ par $\lambda$ est encore de degré $m'$. 
\end{theo}

\subsubsection{} On considère maintenant $(D,F)$ un $\varphi$-module semi-abélien filtré sur $L$ muni d'une donnée de descente pour une extension galoisienne $L/K$ de groupe $G$, c'est-à-dire d'une collection d'isomorphismes $(f_{\sigma})_{\sigma\in G}$ qui vérifie la condition de cocycle $f_{ab}=a(f_b)\circ f_a$ pour $a,b\in G$ et qui sont compatibles à la polarisation $\lambda_0\colon D\to D^t$ et aux suites exactes qui définissent $D$ comme module semi-abélien. Comme $D$ est défini sur $K_0$, cela revient à dire que la filtration $F\subset D_L$ vérifie $f_a(F)=a(F)$ ou encore que $D_L$ est muni d'une action semi-linéaire du groupe opposé $G^{\mathrm{op}}$ de $G$ et que $F$ provient d'un sous-espace de dimension $g$ du $K$-espace vectoriel $D_L^G$ des invariants sous cette action. Par oubli de la filtration admissible $F$ on récupère $D$ et une collection d'automorphismes de celui-ci induisant un morphisme $G\to \Aut_{\varphi} D$. On suppose de plus ce morphisme injectif.

\medskip

D'après le théorème~8.57 de \cite{FO22}, le module $D$ muni de sa donnée de descente correspond à une $\Q_p$-représentation de $G_{K}$, que l'on note $(V, \rho)$. Cette équivalence de catégories est compatible aux suites exactes et à la dualité, donc $V$ est équipé d'un isomorphisme $V\to V^t$ que l'on note $\lambda$ et s'inscrit dans une suite exacte de $\Q_p$-représentations
\[
\begin{tikzcd}
    0\arrow[r] & V_T \arrow[r] & V \arrow[r] & V_B \arrow[r] & 0.
\end{tikzcd}
\]
Par un choix de relèvement à $G_K$ des éléments $h\in G$, que l'on note encore $h$, la donnée de descente $(f_h)_{h\in G}$ correspond aux isomorphismes des $\Q_p$-représentations de $G_L$ donnés par $\rho(h)\colon V\rightarrow V_h$ où $V_h$ est la représentation $(V, \rho(h\cdot h^{-1}))$. On montre en remarque~\ref{rem:indeprelev} que la construction est indépendante du choix de relèvement des éléments de $G$ effectué ici.

On choisit une base convenable $\mB$ adapté à l'inclusion $V_T\subset V$. Soit $\Lambda'$ le réseau engendré par $\mB$. Par continuité de l'action de $G_K$, on obtient un réseau stable $\Lambda$ défini par la somme $\sum_{g\in G_K} g\Lambda'$ tel que $\Lambda_T=\Lambda\cap V_T$ est un réseau stable de $V_T$ et $\Lambda_B=\Lambda/\Lambda_T$ est un réseau stable de $V_B$. Par un choix de base convenable on a $\Lambda=\Lambda_T\oplus \Lambda_B$. Quitte à réduire le réseau $\Lambda_B$ par multiplication par une puissance de $p$ convenable on peut supposer de plus que $\lambda_B(\Lambda_B)\subset \Lambda_B^{vee}$, i.e. $\lambda_B$ définit un morphisme de la $\Z_p$-représentation $T_B$ dans sa duale. On réduit alors la partie $\Lambda_T$ du réseau convenablement pour garder la propriété de stabilité galoisienne. On note encore $\Lambda\subset V$ le réseau stable modifié ainsi et on considère le réseau $\lambda_T(\Lambda_T)\oplus \Lambda_B^{\vee}$ de $V^t$. Comme $\lambda$ commute aux actions galoisiennes c'est un sous-réseau stable de $V^t$ que l'on note $\Lambda^t$. On obtient ainsi deux $\Z_p$-représentations de $G_{K}$ qui s'insèrent dans des suites exactes et qui sont muni d'un morphisme $\lambda\colon \Lambda\to \Lambda^t$ dont la restriction aux parties toriques $\Lambda_T$ et $\Lambda_T^t$ est un isomorphisme. 

\medskip

Le réseau $\Lambda$ muni de sa structure de représentation cristalline de $G_L$, de la même manière que $V$, est équipé d'une donnée de descente $(f_{\sigma})_{\sigma\in G}$ compatible à $\lambda$. L'image $\mfm$ de $\Lambda$ par le foncteur $D_{\mfS}$ dans la catégorie des modules de Breuil-Kisin fournit un groupe $p$-divisible $\Gamma$ sur $\OO_L$ avec une donnée de descente $(f_{\sigma})_{\sigma\in G}$. Par commutativité du diagramme au paragraphe~\ref{diag:categories} et compatibilité au passage aux fibres spéciales on déduit que le module de Dieudonné $M$ de la fibre spéciale $\Gamma_k$ de $\Gamma$, qui est le sous-réseau de $D$ obtenu par le foncteur $\mD$, vérifie $D=M\otimes K_0$ comme $\varphi$-module. Le module de Dieudonné $M$ s'identifie donc à un sous-réseau stable pour $G$ de $D$ et la donnée de descente sur $\Gamma$ induit l'injection $G\hookrightarrow \End M \subset \End_{\varphi} D$. Tout cela vaut de même pour $\Lambda^t$, d'où un groupe $p$-divisible $\Gamma^t$ et un morphisme 
\[
\lambda\colon \Gamma\to \Gamma^t.
\]
Ce morphisme coïncide avec $\lambda_0$ après passage à la fibre spéciale et au foncteur \og{}module de Fontaine\fg{} par construction. 

\begin{rem} \label{rem:indeprelev} Deux choix de relèvements d'un élément $a\in G_K$ diffèrent par un élément $h\in G_L$ et on a un triangle commutatif de morphismes de représentations
\[
\begin{tikzcd}
    &  V^{ah}\\
    V \arrow[r, "\rho(a)"'] \arrow[ur, "\rho(ah)"] & V^a \arrow[u, "\rho(aha^{-1})"'] 
\end{tikzcd}
\]

Or par définition du foncteur de Fontaine $\Rep_{\Q_p}(G_L) \rightarrow \operatorname{MF}^{\varphi}$ les morphismes $\{\rho(h)\}_{h\in G_L}$ ont pour image l'identité dans la catégorie des $\varphi$-modules. Comme cette catégorie est, à isogénie près, celle des fibres spéciales des groupes $p$-divisibles on en déduit que le morphisme induit par $\rho(aha^{-1})$ sur $\Gamma_k$ est l'identité. La fidélité du passage à la fibre spéciale assure alors que l'image du morphisme $\rho(aha^{-1})$ dans la catégorie des groupes $p$-divisibles sur $\OO_L$ est encore l'identité. 
\end{rem}

\begin{rem}
    On peut obtenir $\mfm$ comme image par la composition de deux foncteurs. Le premier est l'inverse du foncteur $\widehat{T}_{L/K}$, donnant une équivalence de catégories d'après le théorème~3.5 \cite{Oze17}, entre la catégorie $\Rep_{\Z_p}^{pst, L,[0,r]}(G_K)$ des réseaux dans les représentations $p$-adiques de $G_K$ potentiellement semi-stables sur $L$ à poids de Hodge-Tate dans $[0,r]$ et la catégorie $\operatorname{Mod}_{\mfS}^{\varphi,r,\widehat{G_L},K}$ des $(\varphi, \widehat{G_K},K)$-modules dont les objets sont des triplets $(\mfm,\varphi, \widehat{G_L})$ où $\mfm$ est un module de Kisin de hauteur plus petite que $r$, voir la définition~3.4 du même article pour les détails. Le deuxième est le foncteur d'oubli $\operatorname{Mod}_{\mfS}^{\varphi,r,\widehat{G_L},K}\rightarrow \operatorname{BT}_{\mfS}^{\varphi}$, la construction de Liu et Ozeki étant compatible aux foncteurs de Kisin. 
\end{rem}

La construction faite à ce paragraphe est résumée par l'énoncé suivant.

\begin{prop}\label{prop:phimodilftrintopdiv} Soient $\overline{k}$ la clôture algébrique d'un corps fini $k$ de caractéristique $p$ et $K_0$ le corps des fractions de l'anneau des vecteurs de Witt $W(\overline{k})$. Soit $K/K_0$ une extension finie et $L/K$ une extension galoisienne de groupe $G$. Soit $D$ un $\varphi$-module semi-abélien sur $K_0$ filtré admissible sur $L$ muni d'une donnée de descente $(f_{\sigma})_{\sigma\in G}$ pour l'extension $L/K$ qui fournit une injection $G\hookrightarrow \Aut_{\varphi} D$. Alors, il existe des groupes $p$-divisibles $\Gamma$ et $\Gamma_t$ sur l'anneau des entiers $\OO_L$ de $L$ qui vérifient les propriétés suivantes.

\begin{enumerate}
    \item Les groupes $p$-divisibles $\Gamma$ et $\Gamma^t$ sont munis de données de descentes pour l'extension $\OO_L/\OO_{K}$. Par restriction à la fibre spéciale et passage au foncteur \og{} module de Fontaine \fg{} la donnée de descente sur $\Gamma$ induit l'injection $G\hookrightarrow \Aut_{\varphi} D$.
    \item Les modules de Dieudonné $M_{\Gamma}$ et $M_{\Gamma^t}$ des fibres spéciales de $\Gamma$ et $\Gamma^t$ vérifient $M_{\Gamma}\otimes_{W(k)} K_0=D$ et $M_{\Gamma^t}\otimes_{W(k)} K_0=D^t$.
    \item Il existe un morphisme 
    \[
    \lambda\colon \Gamma\to \Gamma^t
    \]
    telle que sa restriction aux fibres spéciales coïncide avec $\lambda_0$ par passage au foncteur \og{}module de Fontaine\fg{}. De plus, $\lambda$ commute aux données de descentes sur $\Gamma$ et $\Gamma^t$. 
\end{enumerate}

tel que le module de Dieudonné $M$ de la fibre spéciale de $\Gamma$ vérifie $M\otimes_{W(k)} K_0=D$ et $\Gamma$ est muni d'une donnée de descente $(f_{\sigma})_{\sigma\in G}$ pour l'extension $L/K$ qui induit l'injection  par passage à la fibre spéciale.
\end{prop}

\begin{rem}
    La donné de descente demandée en hypothèse de la proposition~\ref{prop:phimodilftrintopdiv} correspond à la propriété (\ref{enum:stabgal}) des filtrations données par le théorème~\ref{theo:EAdm}. 
\end{rem}

\subsubsection{} On peut finalement démontrer le théorème \ref{theo:deform}.

\begin{proof}
    Soit $D$ le $\varphi$-module issu de $A_0$, qui est semi-abélien polarisé suivant le paragraphe~\ref{subsub:liensemiabphimod}. D'après le théorème~\ref{theo:EAdm}, pour une extension galoisienne $L/K$ de groupe de Galois $G$ et une filtration admissible $F$ sur $D_L$ telle que l'action de $G$ sur le $\varphi$-module filtré $(D,F)$ soit une donnée de descente pour $L/K$. On peut donc appliquer la proposition~\ref{prop:phimodilftrintopdiv} pour obtenir des groupes $p$-divisibles $\Gamma$ et $\Gamma^t$ sur $\OO_L$ qui s'inscrivent dans des suites exactes courtes
    \[
    \begin{tikzcd}
          \Gamma_T \arrow[r] & \Gamma \arrow[r] & \Gamma_B & \mathrm{et} & \Gamma^t_T \arrow[r] & \Gamma^t \arrow[r] & \Gamma^t_B
    \end{tikzcd}
    \]
    et qui sont muni d'un morphisme $\Gamma \to \Gamma^t$ respectant ces suites et induisant un isomorphisme $\Gamma_T\to \Gamma^t_T$. De plus, par construction la fibre spéciale $\Gamma_k$ est isogène à $A_0[p^{\infty}]$ et $\Gamma$ est muni d'une donnée de descente $(f_{\sigma})_{\sigma\in G}$ qui induit une injection $G\hookrightarrow \Aut \Gamma_k \subset \Aut_{\varphi} D$. Le quotient de $A_0$ par le noyau de l'isogénie $A_0[p^{\infty}]\rightarrow \Gamma_k$ est une variété semi-abélienne $C_0$ dont le groupe $p$-divisible est $\Gamma_k$ par construction. De la même façon, on obtient $C_0^t$ et une polarisation $C_0\to C_0^t$. Il suit, d'après le théorème~3.2.1 de \cite{BM19}, que $C_0$ se relève en un schéma semi-abélien formel $\mC$ sur $\OO_L$ muni d'une polarisation $\lambda$ et d'une donnée de descente $(f_{\sigma})_{\sigma\in G}$ compatible à $\lambda$ qui induit l'injection $G\subset \Aut A_0$ sur la fibre spéciale. Par le théorème d'algébrisation de Grothendieck (\cite{EGA3}, théorème~5.4.5) le schéma formel $\mB$, partie abélienne de $\mC$ est algébrisable en un schéma $B$ sur $\OO_L$ muni de la même structure. Il suit, de la même façon qu'au paragraphe~1 du chapitre~2 de \cite{FC90} que $\mC$ est algébrisable en un schéma semi-abélien $A$ sur $\OO_L$ muni à nouveau de la même structure de descente.

    \medskip

    Pour l'assertion sur le degré de la polarisation on remarque que celle-ci s'obtient de la construction comme relevé de la polarisation $C_0\to C_0^t$ déduite des $p$-isogénies entre $A_0$, $C_0$, $A_0^t$ et $C_0^t$. La partie première à $p$ du degré est donc inchangée et la $p$-partie est fournie par l'indice de l'inclusion des réseaux stables $\Lambda\subset \Lambda^t$. L'assertion sur le degré torique s'obtient de la même façon.
\end{proof}

\begin{rem}
    Dans le cas maximalement dégénéré, qui correspond à celui où $A_0$ est un tore déployé sur $k$, la condition d'admissibilité est triviale comme on l'a vu. Il est notable que la construction par déformation présentée dans cette partie se fasse ici directement et sans difficulté, grâce à la rigidité des tores. Le groupe algébrique $A_0$ est un tore déployé sur $k$, donc $A_0\simeq \G_m^n$, le relèvement du groupe $G$ en $K$-automorphismes d'une variété semi-abélienne en caractéristique $0$ se fait simplement en considérant $\G_m^n$ comme schéma en groupes sur $\OO_K$ avec l'inclusion $G\subset \GL_n(\Z)= \Aut \G_m^n$ qui relève l'inclusion donnée sur la fibre spéciale. La donnée de descente voulue est alors celle donnée par l'inclusion $\Gal(L/K)\simeq G \subset \Aut_L \G_m^n$ qui vérifie en particulier la condition de cocycle.
\end{rem}

\subsection{Dégénérescence et monodromie finie} \label{sub:degmonod}

On se replace désormais dans la situation de départ. Soit un groupe fini  $G$, de ramification en $p$, avec une injection $\iota \colon G\hookrightarrow \Aut (A_0, \lambda_0)$ pour une variété semi-abélienne polarisée $A_0$ de dimension $g$ sur la clôture algébrique $\overline{k}$ d'un corps fini $k$ de caractéristique $p$. On dispose d'après le théorème~\ref{theo:deform} d'une extension galoisienne $L/K$ de groupe $G$ de corps résiduel $\overline{k}$ et d'un schéma semi-abélien $C$ sur $L$ muni d'une donnée de descente $(f_{\sigma})_{\sigma\in G}$ qui induit l'injection $\iota$ par passage à la fibre spéciale ainsi que d'un relevé $\lambda$ de $\lambda_0$ compatible aux $(f_{\sigma})_{\sigma\in G}$.  

\smallskip

On se ramène à une variété abélienne $A$ sur $L$ munie d'une donnée de descente par la théorie de Faltings et Chai. Cela correspond à montrer que l'on peut ajouter une donnée de dégénérescence compatible -- au sens du chapitre 3 de \cite{FC90} -- à la donnée de descente sur les schémas semi-abéliens fournis par le théorème~\ref{theo:deform}.

\begin{prop} \label{prop:degenvar} Soit $C$ un schéma semi-abélien sur $\OO_L$ dont la partie torique est déployée muni d'une polarisation $\lambda$ qui est un isomorphisme sur les parties toriques. On suppose de plus $C$ muni d'une donnée de descente $(f_{\sigma})_{\sigma\in G}$ compatible à $\lambda$ qui induit une injection $\iota\colon G\rightarrow \Aut (C_k, \lambda_k)$. Alors, il existe une variété abélienne $A$ sur $L$, dont la réduction $A_k$ est isomorphe à $C_k$, munie d'une donnée de descente $(f_{\sigma})_{\sigma\in G}$ qui induit l'injection $\iota$ ainsi que d'une polarisation de degré $\deg \lambda$.  
    
\end{prop}
\begin{proof}
    Le schéma semi-abélien $C$ avec la polarisation de sa fibre spéciale correspond à une donnée de dégénérescence polarisée, c'est-à-dire un objet de la catégorie $\mathrm{DD}_{\mathrm{pol}}$ définie p.57--58 de \cite{FC90}, à l'exception du morphisme $\underline{Y}\rightarrow C(L)$ du point $(vi)$. Le morphisme cherché $\underline{Y}\rightarrow C(L)$ correspond à un relevé de $c^t\colon \underline{Y}\to B(L)$ qui vérifie une condition de symétrie vis-à-vis de $\lambda$, une condition de positivité et une condition de stabilité pour la donnée de descente. Montrons l'existence d'un tel morphisme, avec les notations de Faltings et Chai. Du fait que la partie torique est déployée, les faisceaux $\underline{X}$ et $\underline{Y}$ s'identifient à des $\Z$-modules libres de rang $n$ que l'on note $X$ et $Y$. On introduit plusieurs ensembles. L'ensemble $R$ des relevés de $c^t$ avec comme sous-ensemble $R_s$ l'ensemble des relevés symétriques. L'hypothèse sur $\lambda$ se traduit ici par le fait que $\Phi\colon Y\to X$ est un isomorphisme. Montrons que cela assure que l'ensemble $R_s$ est non vide. On choisit pour cela des bases $(x_i)_{i\in \{1,\dots, n\}}$ et $(y_i)_{i\in \{1,\dots, n\}}$ de $X$ et $Y$ adaptés au morphisme $\Phi$. Autrement dit, $\Phi(y_i)=x_i$. On choisit de plus un relevé arbitraire $u\in R$, qui correspond à un choix compatible de trivialisations des faisceaux $c^t(y_i)^*\mL_{x_i}$ où $\mL_{x_i}$ est le faisceau inversible sur $B$ correspondant à $c(x_i)\in B^{\vee}$. La structure sur $C$ fourni des isomorphismes de symétrie
    \[
    s_{ij}\colon c^t(y_i)^*\mL_{\Phi(y_j)} \simeq c^t(y_j)^*\mL_{\Phi(y_i)}
    \]
    qui s'identifient par notre choix de relevé à des éléments $s_{ij}\in \OO_L^{\times}$. Un choix symétrique de relevé de $c^t$ correspond alors à un choix de trivialisations $(\alpha_{ij})_{i,j\in \{1,\dots, n\}}$ des fibres génériques des faisceaux $c^t(y_i)^*\mL_{x_j}$, et par identification à une matrice à coefficients dans $L^{\times}$, qui vérifient
    \[
    \alpha_{ij}= \alpha_{ji} \circ s_{ij}.
    \]
    Toujours avec les identifications dû à notre choix arbitraire de relevé et à l'hypothèse sur $\Phi$ cela revient simplement à choisir des éléments $\alpha_{ij}\in L^{\times}$ tels que l'égalité
    \[
    \alpha_{ij}=s_{ij}\alpha_{ji}
    \]
    est vérifiée. Il est clair que cela est possible.
    
    \medskip
    
    Pour $u,v$ deux éléments de $R$ la différence $u-v$ appartient au $\Z$-module libre
    \[
    \Hom(Y, T(L)) = \operatorname{Bil}(Y\times X, L^{\times})= \Hom(Y\otimes X, L^{\times})=Y^{\vee}\otimes X^{\vee}\otimes L^{\times}
    \]
    que l'on note $V$. On note $V_s$ la partie de $V$ formée des éléments symétriques, c'est-à-dire les $\psi\colon Y\times X \to L^{\times}$ tels que $\psi\circ (\mathrm{id}_Y\times \Phi)$ est symétrique au sens usuel. De plus, les ensembles $R$ et $R_s$ sont des espaces homogènes principaux sur $V$ et $V_s$ respectivement.

    \medskip
    
    La donnée de descente sur $C$, que l'on identifie ici a une action de $G$ sur $C$, compatible à son action naturelle sur $\OO_L$, induit des actions sur tous les ensembles $R$, $R_s$, $V$ et $V_s$. Cette action s'écrit explicitement grâce aux actions données sur $X$, $Y$ et $L$. En particulier, pour $u\in R$ et $\sigma\in G$ on défini un élément $\sigma\cdot u -u$ de $V$. On obtient ainsi un cocycle $G\to V$ et changer $u$ par l'ajout d'un élément $v$ de $V$ revient à lui ajouter un cobord. Autrement dit, la donnée de descente détermine un élément $\gamma\in H^1(G,V)$. Il existe un $u\in R$ stable pour cette donnée si et seulement si cette élément s'annule. D'après le corollaire p.170 de \cite{Se62} et avec la description $V=Y^{\vee}\otimes X^{\vee}\otimes L^{\times}$ le groupe $H^1(G,V)$ s'annule. De la même manière dans le cas symétrique, en remarquant que $V_s$ s'identifie à $S^{\vee}\otimes L^{\times}$ où $S=Y\otimes X/\langle y\otimes \Phi(y')-y'\otimes \Phi(y)\rangle$ on obtient l'annulation de $H^1(G,V_s)$. On en déduit l'existence d'un relevé symétrique $u\in R_s$ de $c^t$ compatible à la donnée de descente. Il reste à voir que l'on peut le modifier pour qu'il vérifie la propriété de positivité. 

    \medskip

    La $\OO_L$-structure définit une application $b\colon R\to \operatorname{Bil}(Y\times Y, \R)$, qui s'interprète avec un choix de base $(y_i)_{i\in \{1,\dots, n\}}$ de $Y$ et la donnée d'un relevé symétrique de $c^t$ comme l'association à un élément $u\in R$ de la matrice des valuations des trivialisations données par $u$ avec les choix en vigueur. En particulier, la partie $R_s$ des relevés symétriques a son image dans les applications bilinéaires symétriques. Pour un élément $v\in V$ on obtient un élément $b(v)\in \operatorname{Bil}(Y\times Y, \R)$ par la composition
    \[
    \begin{tikzcd}
        Y\times Y\arrow[r, "\mathrm{id}_Y\times \Phi"] & Y\times X \arrow[r, "v"] & L^{\times} \arrow[r] & \R.
    \end{tikzcd}
    \]
    Il est clair que $b(u+v)=b(u)+b(v)$. Pour $u\in R_s$ stable pour la donnée de descente on construit un élément vérifiant la propriété de positivité en prenant $u+nv$ pour un $n$ assez grand et $v\in V_s^G$ tel que $b(v)$ est défini positif. Un tel élément $v$ s'obtient, par exemple, par la matrice $(p^{m_{ij}})_{i,j\in \{1,\dots, n\}}$ où $(m_{ij})_{i,j\in \{1,\dots, n\}}\in \mathrm{M}_n(\Z)$ est symétrique définie positive et stable par $G$.

    \medskip

    Par l'équivalence de catégories donnée au corollaire~7.2 du chapitre~3 de \cite{FC90}, le schéma $C$ muni de cette structure additionnelle fourni un schéma en groupes commutatif $\mB$ sur $\OO_L$ de fibre générique une variété abélienne polarisée $A$ de dimension $g$, dont la polarisation est de degré $\deg \lambda$. La fonctorialité assure de plus que $\mB$ est aussi muni d'une donnée de descente $(f_{\sigma})_{\sigma\in G}$. L'extension de Raynaud de $\mB$ est $C$ et la fibre spéciale de $\mB$ est la variété semi-abélienne $A_0$. De plus, par restriction à la fibre spéciale, la famille $(f_{\sigma})_{\sigma\in G}$ induit l'injection $\iota$.  

    \medskip

    Soit alors $\mA$ le modèle de Néron de $A$ sur $\OO_L$. Le conjugué $\mA^{\sigma}$ de $\mA$ par $\sigma\in G$ est le modèle de Néron de la variété abélienne $A^{\sigma}$. Par la propriété de Néron les morphismes $(f_{\sigma})_{\sigma\in G}$, d'abord restreint aux fibres génériques, s'étendent à ces modèles. Par ailleurs, par la proposition 3.3 de l'exposé IX de \cite{sga} on a des immersions ouvertes $\mB^{\sigma} \hookrightarrow \mA^{\sigma}$ pour tout $\sigma\in \mathrm{Gal}(L/K)$. En particulier, ces immersions induisent des isomorphismes sur les composantes connexes des neutres des fibres génériques et spéciales. On a donc $A_0=\mB_k=\mA^{\circ}_k=\mA^{\sigma \circ}_k$. De plus, l'unicité de l'extension d'un morphisme de $A^{\sigma}$ à $\mA^{\sigma}$ assure que les morphismes $f_{\sigma}$ entre les modèles de Néron des variétés abéliennes $A^{\sigma}$ induisent l'injection $\iota\colon G\hookrightarrow \Aut \mA^{\circ}_k=\Aut A_0$ par restriction aux fibres spéciales.
\end{proof}

On peut finalement conclure que $G$ est le groupe de monodromie finie de la variété abélienne $B$ déduite de $A$ par la donnée de descente $(f_{\sigma})_{\sigma\in G}$ (dans le cas des courbes elliptiques avec bonne réduction ce résultat est dû à M. Volkov, voir le théorème~4.5 de \cite{Vo01}).

\begin{theo}\label{theo:descentemonod}
Soit $A$ une variété abélienne munie d'une polarisation de degré $d$ sur $L$ munie d'une donnée de descente $(f_{\sigma})_{\sigma\in G}$ relative à l'extension $L/K$. On suppose que
\begin{itemize}
\item[$(i)$] $A_L$ a réduction semi-stable.
\item[$(ii)$] On a une injection $\mathrm{Gal}(L/K)\rightarrow \Aut A_k$ induite par la donnée de descente $(f_{\sigma})_{\sigma\in G}$. 
\end{itemize}
Alors $A$ descend en une variété abélienne $B$ munie d'une polarisation de degré $d$ sur $K$ et $G$ est le groupe de monodromie finie de $B$.
\end{theo}
\begin{proof}
La donnée de descente $(f_{\sigma})_{\sigma\in G}$ est effective et $A$ descend à une variété abélienne $B$ sur $K$ munie d'une polarisation de degré $d$. Soit $I_B$ le groupe de Galois absolu de $K_B$, l'extension de $K$ donnée en définition~\ref{def:monod}. Comme $A$ a réduction semi-stable on a $G_L\subset I_B$ et une surjection $\Gal (L/K) \rightarrow \Phi_B$. Il suit d'après \cite{SZ98} Proposition~4.3 (ii) les égalités $(\operatorname{T}_{\ell} B)^{I_B}=(\operatorname{T}_{\ell} B)^{G_L}=\operatorname{T}_{\ell} A_k$ pour $\ell$ premier différent de $p$. Par ailleurs, l'action naturelle de $G_L$ sur $(\operatorname{T}_{\ell} B)^{G_L}$ se factorise par le quotient $\Gal(L/K)$. Or l'action de ce dernier est déduite de la donnée de descente  $(f_{\sigma})_{\sigma\in G}$. D'après $(ii)$, on a donc un diagramme commutatif
\[\begin{tikzcd}
 G_K \arrow[rr] \arrow[dd] &  & \Aut (\operatorname{T}_{\ell} B)^{I_B} \\
  &    \Aut A_k \arrow[ru, hook] & & \\
 \Gal (L/K)  \arrow[r] \arrow[ru, hook] & \Phi_{B}  \arrow[u, hook]  \arrow[ruu, hook] & 
\end{tikzcd}\] 

et la surjection $\Gal (L/K) \rightarrow \Phi_B$ est un isomorphisme ou encore $I_B=G_L$.

\end{proof}

On termine par le résultat principal de cette partie, aboutissement de la construction présentée aux parties \ref{sub:deformHodge}, \ref{sub:degmonod} et qui démontre le théorème~\ref{theo:main} de l'introduction. 

\begin{theo}\label{theo:resumconstruction} Soit $G$ un groupe fini de ramification tel qu'il existe une variété semi-abélienne polarisée $(A_0,\lambda_0)$ sur un corps fini $k$ avec une inclusion $G\subset \Aut (A_0,\lambda_0)$. Alors il existe une variété abélienne $B$ de dimension $\dim A_0$ sur un corps $p$-adique $K$ telle que le groupe de monodromie finie de $B$ en la place $v$ de $K$ soit $G$ et dont la réduction $(B_L)_{\overline{k}}$ de $B_L$ est géométriquement isogène à $A_0$. De plus, si $\lambda_0$ est de degré $m$ alors $B$ est munie d'une polarisation de degré $mp^n$ pour un $n\geq 0$.
\end{theo}
\begin{proof}
    Quitte à quotienter $A_0$ par le noyau de $\lambda_{T_0}$, on peut supposer que $\lambda_{T_0}$ est un isomorphisme. Alors d'après le théorème~\ref{theo:deform} on obtient un schéma semi-abélien $C$ de fibre spéciale $A_0$ qui vérifie les hypothèses de la proposition~\ref{prop:degenvar} où l'injection $\iota \colon G\rightarrow \Aut (C_k, \lambda_k)$ induit l'inclusion $G\subset \End A_0 \otimes \Q$ issue de l'énoncé. Il suit l'existence d'une variété abélienne polarisée $A$ sur $L$ muni d'une donnée de descente pour l'extension $L/K$ de telle sorte que les morphismes de la donnée de descente induisent une injection $\iota\colon\mathrm{Gal}(L/K)\rightarrow \Aut A_{\overline{k}}$ où $A_{\overline{k}}$ est la réduction de $A$. L'injection $\iota$ correspond à nouveau à l'inclusion de l'énoncé vue dans $\End (A_0)_{\overline{k}}\otimes \Q = \End A_{\overline{k}} \otimes \Q$ avec les isomorphismes $G\simeq \Gal (L/K)$ et $A_{\overline{k}}\simeq C_{\overline{k}}$. On est alors en position d'appliquer le théorème \ref{theo:descentemonod} ce qui fournit une variété abélienne $B$ de dimension $\dim A_0$ sur le corps $K$ avec $G$ pour groupe de monodromie finie. Il est clair que $B$ descend à une extension finie $K'$ de $\Q_p$ avec même groupe de monodromie finie.  
    
\end{proof}

Le passage à un corps de nombres s'effectue par une adaptation du théorème~4.3 de \cite{Ph221} où l'on considère une composante irréductible de $H_{g,d}$ de l'espace des variétés abéliennes ayant une polarisation de degré $d$ et linéairement rigidifiées contenant $B$, où $d$ est le degré d'une polarisation sur $B$. 

\begin{cor}
    Soit $G$ un groupe fini de ramification tel qu'il existe une variété semi-abélienne polarisée $(A_0,\lambda_0)$ sur un corps fini $k$ avec une inclusion $G\subset \Aut (A_0,\lambda_0)$. Alors il existe une variété abélienne $B$ de dimension $\dim A_0$ sur un corps de nombres $K$ et une place non archimédienne $v$ de $K$ de corps résiduel $k$ telle que le groupe de monodromie finie de $B$ en $v$ soit $G$. De plus, si $\lambda_0$ est de degré $m$ alors $B$ est munie d'une polarisation de degré $mp^n$ pour un $n\geq 0$.
\end{cor}

\begin{proof}
    Cela suit directement d'une adaptation du théorème~4.3 de \cite{Ph221} appliquée avec la variété abélienne fournie par le théorème~\ref{theo:resumconstruction}. 
\end{proof}

En suivant précisément la construction par approximation faible on peut obtenir de plus que la réduction $(B_L)_w$ de $B_L$ en $w$ est géométriquement isogène à $A_0$ pour une extension $L/K$ et une place $w\mid v$ de $L$ telle que $B_L$ a réduction semi-stable en $w$. 


\begin{ex} \label{ex:Q8Sgreal}
    On poursuit l'exemple~\ref{ex:Q8Sg}. Soit $E$ une courbe elliptique surpersingulière sur $\F_4$. On peut vérifier que son module de Dieudonné $M_E$ vérifie $M_E\otimes_{W(\F_4)} \Q_4= D_{1/2}$ et que son algèbre d'endomorphismes est l'algèbre de quaternions $\mathbf{H}$. L'injection naturelle $\End E^g \hookrightarrow \End D^g$ est donc ici un isomorphisme pour tout $g\geq 1$. De plus, $E$ admet une polarisation principale qui induit $\lambda_0$ sur $D_{1/2}$ et la polarisation produit sur $E^g$ induit $\lambda_0^g$. 

    \smallskip
    
    D'après ce que l'on a vu, on peut appliquer le théorème~\ref{theo:resumconstruction} avec comme donnée de départ l'inclusion $G_g\subset \Aut (A_0,\lambda_0)$ pour $g\geq 1$ ce qui fournit une variété abélienne $A$ sur un corps de nombres avec un groupe de monodromie finie d'ordre $2^{r(2g,2)}$ munie d'une polarisation de degré $2^n$ pour un $n\geq 1$ du fait que $\lambda_0$ est principale. 
\end{ex}

\section{Application au degré de semi-stabilité} \label{sec:degsemistab}

\subsection{Le degré de semi-stabilité} \label{sub:maxd}

\subsubsection{} Le théorème~4.3 de \cite{Ph221} donne une forme de principe local-global pour les groupes de monodromie finie sous l'hypothèse technique que les variétés abéliennes associées doivent être munies d'une polarisation principale. Ce résultat avec la formule, du même article,
\[
d(A)= \underset{v\in \Sigma_K}{\ppcm} \Card \Phi_{A,v}
\]
pour une variété abélienne sur un corps de nombres, ramène la question de la construction de variétés abéliennes avec degré de semi-stabilité maximal à son analogue local. Précisément, si l'on construit des variétés abéliennes principalement polarisées et de même dimension $A_1,\dots, A_n$ sur des corps $p$-adiques $K_1,\dots, K_n$ on obtient l'existence d'une variété abélienne $A$ sur un corps de nombres $L$ avec
\[
\underset{i\in \{1,\dots, n\}}{\ppcm} \Card \Phi_{A_i} \mid d(A). 
\]

Pour un entier $g\geq 1$, la construction CM faite dans \cite{Phi222} - voir la remarque~4.6 - donne que pour tout corps de nombres $K$ il existe une extension finie $L$ de $K$ non ramifiée aux places au-dessus de $2$, des places $v_1,\dots, v_n$ de $L$ au-dessus des diviseurs premiers impairs $p_1,\dots, p_n$ de $M(2g)$, et des variétés abéliennes principalement polarisées $A_1, \dots, A_n$ vérifiant
\[
\Card \Phi_{A_i,v_i} = p_i^{r(2g,p_i)}.
\]

D'après ce qui précède, il est suffisant de construire une variété abélienne principalement polarisée, de dimension $g$ sur un corps $2$-adique avec $\Card \Phi_A=2^{r(2g,2)}$ pour établir l'égalité $d_g=M(2g)$. 

\subsubsection{} La variété abélienne $A$ de l'exemple~\ref{ex:Q8Sgreal} ne convient pas a priori du fait que la polarisation $\lambda$ sur $A$ fournie par notre construction est de degré une puissance $2$. Pour pallier ce défaut on donne un lemme qui permet, quitte à faire une extension finie modérément ramifiée, d'obtenir une variété abélienne isogène à $A$ qui convient. L'invariance par isogénie des groupes de monodromie finie se déduit, par exemple, du fait qu'une isogénie induit un isomorphisme au niveau des modules de Tate sur $\Q_{\ell}$. 

\begin{lem} Soit $\{1\}\neq G$ un groupe fini étale abélien d'ordre une puissance d'un nombre premier $p$ sur un corps $p$-adique $K$. Alors il existe un sous-groupe cyclique $H\subset G$ d'ordre $p$ défini sur une extension finie modérément ramifiée $L$ de $K$.
\end{lem}
\begin{proof}
Il suffit de montrer qu'il existe un tel sous-groupe $H$ fixé par l'action du groupe de ramification sauvage $V_K\subset G_K$. Ce groupe étant un pro-$p$-groupe libre de rang infini, son action sur $G$ est donnée par l'action d'un quotient fini $V$ d'ordre $p^c$ pour un $c\geq 0$. On considère alors l'action induite par $V$ sur l'ensemble $\mC$ des sous-groupes cycliques d'ordre $p$ de $G$. Par un calcul élémentaire, on a 
\[
\Card R= p\Card \mC -\Card \mC +1
\]
où $R$ est le sous-groupe de $G$ engendré par les éléments de $G$ d'ordre $p$. Ce sous-groupe est nécessairement non trivial et d'ordre une puissance de $p$. Il suit que $\Card \mC$ n'est pas divisible par $p$. Comme $\mC$ est l'union disjointe des orbites de ses éléments sous l'action de $V$ il suit qu'au moins une orbite est de cardinal $1$, les orbites non triviales étant de cardinal divisible par $p$. Autrement dit, il existe un sous-groupe cyclique $H\subset G$ fixé par l'action de $V$.   
\end{proof}

\begin{cor}\label{cor:degpoltame}
    Soit $A$ une variété abélienne sur un corps $p$-adique $K$ muni d'une polarisation $\lambda$ de degré $d$ une puissance de $p$. Alors il existe une extension finie modérément ramifiée $L/K$ et une variété abélienne $B$ sur $L$ isogène à $A_L$ telle que $B$ admet une polarisation principale.
\end{cor}

\begin{proof}
    Le résultat s'obtient grâce au lemme par une récurrence sur le degré $d$ de la polarisation. Si $d\geq p$ alors il existe une extension finie $L_1/K$ modérément ramifiée telle qu'un sous-groupe cyclique $H$ d'ordre $p$ du noyau de $\lambda$ soit défini sur $L_1$. La variété quotient $A_{L_1}/H$ est isogène à $A_{L_1}$ et est munie d'une polarisation de degré $d/p$.
\end{proof}

On peut désormais démontrer l'égalité attendue.

\begin{theo}
Pour $g\in \N\setminus\{0\}$ on a 
\[
d_g=M(2g).
\]
\end{theo}
\begin{proof}
D'après le théorème~4.3 de \cite{Ph221} et la remarque~4.6 de \cite{Phi222}, il suffit de construire une variété abélienne $A$ principalement polarisée de dimension $g$ sur un corps $2$-adique avec $\Card \Phi_A=2^{r(2g,2)}$.  

On considère la variété abélienne $A$ sur un corps $2$-adique $K$ obtenue dans l'exemple~\ref{ex:Q8Sgreal}. Le corollaire~\ref{cor:degpoltame} appliqué à $A$ fournit une variété abélienne $B$ principalement polarisée sur une extension $L$ de $K$ finie et modérément ramifiée sur $K$ avec $B$ isogène à $A_L$. Les groupes de monodromie finie étant invariant par isogénie et leur $p$-partie invariante par extension modérément ramifiée, la variété abélienne $B$ convient.  
\end{proof}

\begin{rem}
    On obtient de plus que $d_g$ est aussi le plus petit commun multiple $d_g^{div}$ des entiers $d(A)$ lorsque $A$ varie et non pas seulement le maximum. On en donne rapidement la démonstration. De la formule
    \[
     d(A) = \underset{v\in \Sigma_K}{\ppcm}~\Card \Phi_{A,v}
    \]
    on déduit le fait que $d_g^{div}$ est le plus petit commun multiple des cardinaux des groupes de monodromie finie réalisable en dimension $g$ et l'inégalité $d_g^{div}\leq M(2g)$. Par ailleurs, l'inégalité $d_g\leq d_g^{div}$ est claire. 
\end{rem}

\subsection{Réduction semi-stable déployée} \label{sub:deploy}

On montre ici l'égalité $d_g=d^{d\acute{e}p}$ pour tout entier naturel $g\geq 1$ annoncée dans l'introduction. Cette partie est indépendante des autres et repose sur les propriétés connues des groupes de monodromie finie, en particulier obtenues dans \cite{SZ98}.

\bigskip

On considère un corps de nombres $K$ et une variété abélienne $A$ de dimension $g$ sur $K$. On dit que $A$ a réduction semi-stable déployée sur $K$ si pour toute place $v\in \Sigma_K$ la réduction $A_v$ de $A$ est l'extension d'une variété abélienne par un tore déployé. Il découle facilement du théorème de réduction semi-stable qu'il existe une extension finie $L/K$ telle que $A_L$ a réduction semi-stable déployée. 

\begin{defin}
    On pose 
    \[
    d^{d\acute{e}p}(A)= \min \{ [L:K] \mid A_L \text{ a réduction semi-stable déployé} \}
    \]
    et $d^{d\acute{e}p}_g$ le plus petit entier naturel tel que $d^{d\acute{e}p}(A)\leq d^{d\acute{e}p}_g$ pour toute variété abélienne $A$ de dimension $g$ sur un corps de nombres.
\end{defin}

Il est clair que $d(A)\leq d^{d\acute{e}p}(A)$. Soit $L/K$ une extension finie telle que $A_L$ a réduction semi-stable. Pour une place $w\in \Sigma_L$ le rang torique de la réduction de $A_L$ en $w$ ne dépend que de la place $v$ de $K$ en-dessous de $w$ et on le note $t_v$. Pour déployer un tore de dimension $n$ il suffit d'une extension de degré au plus $M(n)$. On déduit directement la majoration suivante du théorème~2.4 de \cite{Ph221}. 

\begin{prop}\label{prop:formulesc}
    On a 
    \[
    d^{d\acute{e}p}(A) \leq \underset{v\in \Sigma_K}{\ppcm}~ M(t_v) \Card \Phi_{A,v}.
    \]
\end{prop}

\begin{lem}\label{lem:scdivisM}
    Pour tout $g\in \N\setminus\{0\}$ et tout $n\in \{0,\dots, g\}$ on a la divisibilité 
    \[
    M(n)^2M(2g-2n)\mid M(2g).
    \]
\end{lem}
\begin{proof} Cela suit directement de la divisibilité $M(a)M(b)\mid M(a+b)$ pour tous entiers naturels $a$, $b$ qui s'obtient par l'inclusion diagonale $\GL_a(\Q)\times \GL_b(\Q)\subset \GL_{a+b}(\Q).$
\end{proof}

\begin{theo}
    Pour tout $g\in \N\setminus\{0\}$ on a $d^{d\acute{e}p}_g=d_g$.
\end{theo}

\begin{proof}
    L'inégalité $d_g\leq d^{d\acute{e}p}_g$ est triviale. Pour l'autre direction, on considère une variété abélienne $A$ de dimension $g$ sur un corps de nombres $K$. D'après la proposition~\ref{prop:formulesc} on a 
     \[
    d^{d\acute{e}p}(A) \leq \underset{v\in \Sigma_K}{\ppcm}~M(t_v) \Card \Phi_{A,v}.
    \]
    Or d'après \cite{SZ98} corollaire~6.3 le cardinal $\Card \Phi_{A,v}$ divise $M(t_v)M(2g-2t_v)$. Avec le lemme~\ref{lem:scdivisM} il suit $M(t_v) \Card \Phi_{A,v}\mid M(2g)$ d'où l'inégalité
    \[
     d^{d\acute{e}p}(A) \leq M(2g).
    \]
\end{proof}

\printbibliography

\end{document}